\documentclass[11pt,a4paper]{article}
 
\usepackage[margin=2.5cm]{geometry}
\usepackage[T1]{fontenc}
\usepackage[utf8]{inputenc}
\usepackage{lmodern}
\usepackage{setspace}
\usepackage{authblk}

\usepackage[dvipsnames]{xcolor}
\usepackage{graphicx}
\usepackage{tikz}
\usetikzlibrary{math}
\usepackage{pgfplots}
\pgfplotsset{compat=1.18}
 
\usepackage[round]{natbib}
 
\usepackage[figuresright]{rotating}
\usepackage{subcaption}
\usepackage{booktabs}
\usepackage{multirow}
\usepackage{multicol}
\usepackage{placeins}
 
\usepackage{amsmath,amssymb,amsthm}

\theoremstyle{definition}
\newtheorem{example}{Example}

\theoremstyle{remark}

\usepackage{algorithm}
\usepackage{algorithmicx}
\usepackage{algpseudocode}
 
\usepackage{eurosym}
\usepackage{verbatim}
\usepackage{url}
\usepackage{xargs}
\usepackage[colorinlistoftodos,prependcaption,textsize=tiny]{todonotes}
\newcommandx{\GeneralTodo}[2][1=]{\todo[inline,size=\small,linecolor=yellow,backgroundcolor=yellow!50,bordercolor=yellow,#1]{#2}}
\newcommandx{\unsure}[2][1=]{\todo[linecolor=orange,backgroundcolor=orange!25,bordercolor=orange,#1]{#2}}
\newcommandx{\unsureG}[2][1=]{\todo[linecolor=red,backgroundcolor=red!25,bordercolor=red,#1]{#2}}
\newcommandx{\unsureM}[2][1=]{\todo[linecolor=green,backgroundcolor=green!25,bordercolor=green,#1]{#2}}
\newcommandx{\change}[2][1=]{\todo[linecolor=red,backgroundcolor=red!25,bordercolor=red,#1]{#2}}
\newcommandx{\info}[2][1=]{\todo[linecolor=blue,backgroundcolor=blue!25,bordercolor=blue,#1]{#2}}
\newcommandx{\improve}[2][1=]{\todo[linecolor=Plum,backgroundcolor=Plum!25,bordercolor=Plum,#1]{#2}}
\newcommandx{\answer}[2][1=]{\todo[linecolor=orange,backgroundcolor=white!25,bordercolor=orange,#1]{#2}}

\newcommand{\rev}[1]{{\color{black} #1}}
 
\newcommand{\problemname}{$\textsc{PMSTVP}$}
\newcommand{\milpmodel}{$\textsc{MILP}$}
\newcommand{\matheur}{$\textsc{ILS}$}

\newcommand{\tmax}{{t_{\textrm{max}}}}
\newcommand{\pmax}{{p_{\textrm{max}}}}
 
\newcommand{\keywords}[1]{\par\medskip\noindent\textbf{Keywords:} #1}
 
\usepackage[colorlinks=true,linkcolor=blue,citecolor=blue,urlcolor=blue]{hyperref}
 
\begin{document}
 
\title{On parallel machine scheduling with variable energy consumption functions}
 
\author[a]{Mirko Mucciarini}
\author[b,*]{Giulia Caselli}
\author[b]{Daniele De Santis}
\author[b]{Manuel Iori}
\author[c,d]{Juan Jos\'e Miranda-Bront}
 
\affil[a]{DIEF, University of Modena and Reggio Emilia, Modena 41125, Italy}
\affil[b]{DISMI, University of Modena and Reggio Emilia, Reggio Emilia 42122, Italy}
\affil[c]{Universidad Torcuato Di Tella, Escuela de Negocios, Buenos Aires, Argentina}
\affil[d]{Consejo Nacional de Investigaciones Cient\'ificas y T\'ecnicas, Buenos Aires, Argentina}
 
\date{}
 
\maketitle
 
{\renewcommand{\thefootnote}{}%
\footnotetext{$^{*}$Corresponding author: giulia.caselli@unimore.it.\\
E-mail addresses: mirko.mucciarini@unimore.it (M. Mucciarini), giulia.caselli@unimore.it (G. Caselli), 240265@studenti.unimore.it (D. De Santis), manuel.iori@unimore.it (M. Iori), jmiranda@utdt.edu (J.J. Miranda-Bront).}}

\begin{abstract}
Manufacturing facilities face increasing challenges due to renewable energy integration and regulated pricing, motivating energy-efficient scheduling strategies. This paper studies a scheduling problem in which jobs with time-varying energy consumption must be assigned to parallel machines over a discrete planning horizon, subject to per-period energy limits, with the objective of minimizing total energy costs. A time-indexed mixed-integer linear programming (MILP) formulation that accommodates variable energy consumption profiles is proposed, thereby extending commonly used constant-consumption models. To address larger instances, a matheuristic based on Iterated Local Search is developed, implementing a Variable Neighborhood Descent combining adaptations of classical local search operators and leveraging the MILP formulation to explore large neighborhoods. Computational experiments demonstrate that the proposed approach consistently produces solutions of high quality under variable energy consumption profiles.
\end{abstract}


\maketitle

\section{Introduction}
\label{sec:intro}

According to the World Energy Outlook 2023 \citep{IEA2023}, the industrial sector accounts for around 40\% of global final energy consumption and is the largest emitting sector, responsible for approximately 25\% of direct  CO\textsubscript{2} emissions from the global energy system.
These figures underscore the urgent need to rethink energy supply policies and production models.
In the manufacturing sector especially, optimizing production processes requires the implementation of effective scheduling strategies to meet customer demand while limiting resource consumption.

Motivated by the growing concern for climate change, the scientific community has shown an increasing interest in \textit{energy-efficient} (a.k.a. \textit{energy-aware}) \textit{scheduling} (see, e.g., \citealt{FANG2011234}). This field of study focuses on reducing energy consumption, CO\textsubscript{2} emissions, and overall carbon footprint, among other objectives arising in the context of machine scheduling and flexible manufacturing. 
From a practical perspective, companies view Total Energy Cost (TEC) as a key indicator to evaluate the effectiveness of a production schedule. 
At the same time, energy suppliers attempt to reduce the economic costs and environmental impact associated with the generation of high-energy loads in short periods by flattening demand peaks through pricing policies. One of such policies is the Time-of-Use (TOU) tariff, aimed at introducing incentives to modify users' habits by penalizing peak-hour consumption with higher prices \citep{GAHM2016}.
Public institutions have made energy efficiency a strategic priority, as evidenced by the European Union mandate to achieve net-zero emissions for vehicles by 2035, which introduces additional challenges in the automotive production sector \citep{EU2019GreenDeal, EU2021FitFor55}.

In this paper, we address an energy-efficient scheduling problem motivated by Reinova S.p.a., an Italian company operating in the automotive sector, specialized in the testing of batteries and electronic components for electric vehicles. Before installation into vehicles, these components must undergo rigorous tests to meet international quality and safety standards. 
These tests are performed at Reinova by using sophisticated machinery and include the exposure of batteries to physical stresses and different usage scenarios, including repeated charging and discharging cycles, rapid heating and cooling, vibration, and extreme weather conditions, among others, that typically induce heterogeneous, non-constant energy consumption \citep{haspl2022reduction}.

The company owns a number of testing machines that operate in parallel.
Each test can be seen as a job that has to be scheduled on a given machine, having a fixed duration and known total energy consumption. 
However, opposed to the classical scheduling setup, the energy consumption is not uniform over time; rather, it is variable and dictated by the specific nature of each test.
The total amount of energy (i.e., electric power) used at any given time for all the tests must remain within a given threshold, and this becomes a relevant constraint for the scheduling of the tests.
This energy threshold reflects not only physical infrastructure limits, such as the capacity of the internal electrical distribution network, but also contractual agreements with the energy provider, which may impose penalties for exceeding predefined power usage levels.
Furthermore, the energy required for testing can be sourced through multiple channels. In addition to conventional supply from the electricity grid, the company also generates electricity internally through solar panels installed on the roof of the production facility.
The availability of this self-generated energy is inherently variable, as it depends on the hour of the day, as well as on seasonal factors and weather conditions. These factors can cause substantial fluctuations in energy output even within short time periods.
In addition, self-generated energy cannot be stored; it must either be consumed immediately on-site or, if in excess, sold back to the grid.
The energy threshold constraint applies exclusively to the energy actively consumed by the tests and does not account for excess self-generated energy sold back to the grid.
Energy utilization depends directly on the scheduling of the tests, making it essential to align test operations with energy production to optimize resource use.
Crucially, the interaction between energy demand, scheduling, and the variability of both sources and energy prices introduces financial implications. For this reason, aligning test operations with periods of higher self-generation or lower grid prices can significantly reduce operational costs or even generate profit.

Figure \ref{fig:reinova} depicts a graphical representation of the testing system, involving two machines, and its corresponding electric power consumption over a 24-hour period.
Figure \ref{fig:sub-reinovaA} illustrates the scheduled operation of two identical machines: machine $m_1$ is performing a \textit{battery test}, which involves short, repeated charging cycles of a battery (each characterized by a power peak followed by an interval of low activity), resulting in a highly variable power profile; machine $m_2$ is performing a \textit{climatic test}, testing the maintenance of the battery temperature within a stable range (smaller power fluctuations reflect the periodic adjustments in temperature control). 
\begin{figure}
    \centering
    \begin{subfigure}[t]{\textwidth}
    \centering
    \begin{tikzpicture}[scale=1.0]
\draw (0,0) rectangle (12,2); 

\foreach \x in {0.5,1,1.5,...,12} {
    \draw[gray!30, very thin] (\x,0) -- (\x,2);
}

\draw (0,1) -- (12,1);

\node[left] at (-0.05,0.5) {\footnotesize \(m_1\)};
\node[left] at (-0.05,1.5) {\footnotesize \(m_2\)};

\fill[fill=gray!30] (0,0) rectangle (0.25,0.125);
\fill[fill=gray!30] (0.25,0) rectangle (0.5,0.75);
\fill[fill=gray!30] (0.5,0) rectangle (0.75,0.125);
\fill[fill=gray!30] (0.75,0) rectangle (1,0.75);
\fill[fill=gray!30] (1,0) rectangle (1.25,0.125);
\fill[fill=gray!30] (1.25,0) rectangle (1.5,0.75);
\fill[fill=gray!30] (1.5,0) rectangle (1.75,0.125);
\fill[fill=gray!30] (1.75,0) rectangle (2,0.75);
\fill[fill=gray!30] (2,0) rectangle (2.25,0.125);
\fill[fill=gray!30] (2.25,0) rectangle (2.5,0.75);
\fill[fill=gray!30] (2.5,0) rectangle (2.75,0.125);
\fill[fill=gray!30] (2.75,0) rectangle (3,0.75);
\fill[fill=gray!30] (3,0) rectangle (3.25,0.125);
\fill[fill=gray!30] (3.25,0) rectangle (3.5,0.75);
\fill[fill=gray!30] (3.5,0) rectangle (3.75,0.125);
\fill[fill=gray!30] (3.75,0) rectangle (4,0.75);
\fill[fill=gray!30] (4,0) rectangle (4.25,0.125);
\fill[fill=gray!30] (4.25,0) rectangle (4.5,0.75);
\fill[fill=gray!30] (4.5,0) rectangle (4.75,0.125);
\fill[fill=gray!30] (4.75,0) rectangle (5,0.75);
\fill[fill=gray!30] (5,0) rectangle (5.25,0.125);
\fill[fill=gray!30] (5.25,0) rectangle (5.5,0.75);
\fill[fill=gray!30] (5.5,0) rectangle (5.75,0.125);
\fill[fill=gray!30] (5.75,0) rectangle (6,0.75);
\fill[fill=gray!30] (6,0) rectangle (6.25,0.125);
\fill[fill=gray!30] (6.25,0) rectangle (6.5,0.75);
\fill[fill=gray!30] (6.5,0) rectangle (6.75,0.125);
\fill[fill=gray!30] (6.75,0) rectangle (7,0.75);
\fill[fill=gray!30] (7,0) rectangle (7.25,0.125);
\fill[fill=gray!30] (7.25,0) rectangle (7.5,0.75);
\fill[fill=gray!30] (7.5,0) rectangle (7.75,0.125);
\fill[fill=gray!30] (7.75,0) rectangle (8,0.75);

\fill[fill=gray!50] (0,1) rectangle (0.75,1.5);
\fill[fill=gray!50] (0.75,1) rectangle (1.5,1.625);
\fill[fill=gray!50] (1.5,1) rectangle (2.25,1.5);
\fill[fill=gray!50] (2.25,1) rectangle (3,1.625);
\fill[fill=gray!50] (3,1) rectangle (3.75,1.5);
\fill[fill=gray!50] (3.75,1) rectangle (4.5,1.625);
\fill[fill=gray!50] (4.5,1) rectangle (5.25,1.5);

        \foreach \x/\name in {0/0:00, 2/4:00, 4/8:00, 6/12:00, 8/16:00, 10/20:00, 12/24:00} 
            \node[below, scale=0.7] at (\x,0) {\footnotesize \name}; 
    \end{tikzpicture}
    \caption{Test occupation on machines.}
    \label{fig:sub-reinovaA}
\end{subfigure}
    %
    \begin{subfigure}[t]{\textwidth}
            \centering
            \begin{tikzpicture}[scale=1.0]

    \draw[very thin, gray!15] (0,0) grid[xstep=0.5, ystep=0.5] (12,2);

    \draw[->] (0,0) -- (12,0) node[right, yshift=1pt] {\footnotesize T};
    \draw[->] (0,0) -- (0,2) node[above, xshift=-7.75pt, yshift=0.5pt] {\footnotesize kWh};

        \foreach \x/\name in {0/0:00, 2/4:00, 4/8:00, 6/12:00, 8/16:00, 10/20:00, 12/24:00} 
            \node[below, scale=0.7] at (\x,0-0.05) {\footnotesize \name};

    \foreach \y/\label in {0.5/50, 1/100, 1.5/150, 2/200}
        \draw (-0.05,\y) node[left, scale=0.7] {\footnotesize \label};

        \draw[dashed] (0,1.5) -- (12,1.5);

    \draw[thick, black] 
        (0,0.625) -- (0.25,0.625) 
        -- (0.25,1.25) -- (0.5,1.25) 
        -- (0.5,0.625) -- (0.75,0.625) 
        -- (0.75,1.375) -- (1,1.375) 
        -- (1,0.75) -- (1.25,0.75) 
        -- (1.25,1.375) -- (1.5,1.375) 
        -- (1.5,0.625) -- (1.75,0.625) 
        -- (1.75,1.25) -- (2,1.25) 
        -- (2,0.625) -- (2.25,0.625) 
        -- (2.25,1.375) -- (2.5,1.375) 
        -- (2.5,0.75) -- (2.75,0.75) 
        -- (2.75,1.375) -- (3,1.375) 
        -- (3,0.625) -- (3.25,0.625) 
        -- (3.25,1.25) -- (3.5,1.25) 
        -- (3.5,0.625) -- (3.75,0.625) 
        -- (3.75,1.375) -- (4,1.375) 
        -- (4,0.75) -- (4.25,0.75) 
        -- (4.25,1.375) -- (4.5,1.375) 
        -- (4.5,0.625) -- (4.75,0.625) 
        -- (4.75,1.25) -- (5,1.25) 
        -- (5,0.625) -- (5.25,0.625) 
        -- (5.25,0.75) -- (5.5,0.75)
        -- (5.5,0.125) -- (5.75,0.125)
        -- (5.75,0.75) -- (6,0.75)
        -- (6,0.125) -- (6.25,0.125)
        -- (6.25,0.75) -- (6.5,0.75)
        -- (6.5,0.125) -- (6.75,0.125)
        -- (6.75,0.75) -- (7,0.75)
        -- (7,0.125) -- (7.25,0.125)
        -- (7.25,0.75) -- (7.5,0.75)
        -- (7.5,0.125) -- (7.75,0.125)
        -- (7.75,0.75) -- (8,0.75)
        -- (8,0) -- (12,0)
        ;

        \begin{scope}[shift={(9,1)}] 
            \draw[black, thick] (0,0) -- ++(1,0) node[right] {\scriptsize {Total consumption}};
            \draw[dashed] (0,-.3) -- ++(1,0) node[right] {\scriptsize {Energy budget}};
        \end{scope}
        
            \end{tikzpicture}
            \caption{Electric power consumption.}
            \label{fig:sub-reinovaB}
        \end{subfigure}
    %
    \caption{Graphical representation of the testing system with independent test occupation on two machines (Fig.~\ref{fig:sub-reinovaA}), and the overall electric power consumption and energy budget (Fig.~\ref{fig:sub-reinovaB}).}
    \label{fig:reinova}
\end{figure}
In Figure \ref{fig:sub-reinovaB}, the threshold is represented by the dashed horizontal line, and the total power consumption curve of the facility is the black step line below, aggregating the energy used by both machines at each time period.
This simplified representation of the Reinova testing environment demonstrates how the variability in test energy consumption can impact the energy management strategy of the facility.

This application introduces several interesting characteristics, not necessarily restricted to the specific case study but widely applicable to other contexts related to energy-efficient scheduling and, more generally, to scheduling under renewable resources. 
Variable resource consumption, sourcing, and costs arise, for instance, in data centres optimizing job execution to align with dynamic electricity pricing and renewable availability, in manufacturing plants integrating solar or wind power into their operations, and in smart grid environments where residential or commercial systems must be scheduled to match intermittent energy supply.
Formally, we study the problem of assigning a set of jobs to a set of identical parallel machines, ensuring that each job is entirely processed within a discrete finite time horizon. This time horizon corresponds to a single working day and accounts for variable energy costs, including the costs of buying and the prices of selling energy from the grid per time slot. The objective is to minimize the \rev{TEC}, calculated as the total cost of the energy purchased from the grid minus the total profit from the energy sold to the grid, while respecting energy limit constraints and considering the time-dependent energy consumption functions of jobs, each characterized by its specific, irregular profile. We refer to the problem as Parallel Machine Scheduling with Time-of-Use costs and Variable consumptions Problem (\problemname). To the best of our knowledge, the \problemname\ has not been addressed in the literature.

In the recent survey by \citealt{CATANZARO2023} on energy-aware job scheduling on machines, no studies address non-constant energy consumption profiles for jobs.
Despite the growing body of modern research incorporating time-dependent parameters in machine scheduling problems, especially time-dependent processing times (see \citealt{GAWIEJNOWICZ2025}), the assumption of constant energy consumption during the entire job execution (which we refer to as \textit{rectangular}) remains prevalent.
A recent example is provided by \cite{GAGGERO2023}, who study the bi-objective parallel machine scheduling problem under a TOU tariff policy (BPMSTP), aiming to minimize both makespan and TEC.  The authors propose a compact Mixed Integer Linear Programming (MILP) formulation based on a key property of the BPMSTP: if two jobs have the same, fixed  processing time, exchanging them in a schedule does not affect the TEC or the makespan. For the \problemname, 
this property does not hold over the set of feasible schedules because variable consumption can not only affect the TEC but, moreover, such an exchange may be infeasible due to the maximum energy constraint. 
Methodologically, considering variable consumption complicates the design of heuristics, as more complex validations are required to ensure feasibility, as well as the \milpmodel\ models, since formulations need to account for the consumption over time explicitly. 
Therefore, new modeling frameworks are required to address the additional difficulty introduced by variable consumption patterns (i.e., variable resource consumptions of resources during the processing of jobs), not considered in the existing ones based on \textit{rectangular} jobs, and alternative solution methods are needed. 
At the same time, this added difficulty reflects real-world challenges, as demonstrated by the Reinova case study and other industrial applications, motivating the study of the \problemname\ within the machine scheduling research field.

Our primary objective is to fill this gap in the literature by developing mathematical models and heuristic algorithms for the PMSTVP, and more in general for energy-efficient scheduling under time-dependent costs and variable consumptions, and to provide managerial insights regarding the impact of these characteristics into the scheduling process.

The main contributions of this paper are as follows:
\begin{itemize} 
    \item We introduce the \problemname\ and propose a general modelling framework for parallel-machine scheduling with time-dependent resource consumption under variable energy sourcing and pricing;
    \item We propose a MILP formulation and develop a tailored matheuristic capable of generating high-quality solutions for relatively large instances\rev{, together with a properly calibrated Simulated Annealing (SA) metaheuristic used as a benchmark to assess the matheuristic's performance};
    \item We conduct extensive computational experiments to evaluate the effectiveness of the proposed methods. Additionally, we design a tailored experimental framework to analyze the impact of explicitly accounting for variable consumption functions within our scheduling problem, and discuss managerial implications.
\end{itemize}

The remainder of the paper is organized as follows. In Section~\ref{sec:lit}, we provide a summary of the scheduling literature most related to the \problemname.
Section~\ref{sec:problem} defines the \problemname\ and introduces the notation used throughout the paper. Section~\ref{sec:model} presents the MILP formulation, while Section~\ref{sec:heur} details the matheuristic. Section~\ref{sec:exp} describes the experimental setup and summarizes the computational results. Finally, Section~\ref{sec:concl} concludes the paper and outlines potential directions for future research.

\section{Literature review}
\label{sec:lit}


Scheduling problems are among the most studied problems in the Operations Research (OR) literature, with countless variants covering a wide spectrum of practical applications arising in manufacturing, production, and service industries. 
We refer to \cite{pinedo2022scheduling} for an updated review of theoretical models and practical applications, and to \cite{agnetis2014models}, \cite{strusevich2017scheduling}, and \cite{gawiejnowicz2020models} for recent books on modern scheduling.
In general, scheduling problems are $\mathcal{NP}$-hard and solution algorithms usually resort to mathematical programming and heuristics (see, e.g., \citealt{ying2024revisiting}). Due to the extensive literature, in what follows we concentrate on recent contributions on energy-efficient and resource-constrained machine scheduling problems, and (most related) time-dependent scheduling problems \citep{gawiejnowicz2020models}.

\subsection{Energy-efficient parallel machines} 

\cite{GAHM2016} provide an extensive literature review on Energy-Efficient Scheduling (EES), covering both methodological and empirical studies on the impact of energy-aware practices among manufacturing companies, and propose a framework to classify a growing body of interdisciplinary literature. More recently, \cite{CATANZARO2023} survey EES approaches under TOU-based tariff policies, developing a taxonomy for the related literature and discussing both mathematical models and alternative solution methodologies.

Energy-saving policies can be implemented in different ways. One alternative is to avoid using all the available machines at the same time and instead restrict the schedule to be executed using only a subset of them at any given time. This strategy is commonly referred to as the not-all-machine policy \citep{FengPeng2024}. Another option is to limit energy consumption by explicitly limiting the feasibility of a schedule through tailored energy-related constraints. 
\cite{LiEtAl2024} propose an SA for a scheduling problem on parallel machines with machine-dependent energy consumption costs where the objective is to minimize the makespan under a total energy consumption cost constraint. 
A related work is proposed by \cite{TIAN2024}, who consider a single machine parallel batch scheduling problem under TOU tariffs. The objective is to minimize the TEC while an upper limit on the makespan is imposed as a feasibility constraint. The authors propose a branch-and-price approach and a column-generation-based heuristic.
\cite{li2024energy} study the problem of identical parallel machine scheduling under peak power consumption and deadline constraints, where the objective is to maximize the total value of the selected jobs. The authors show that the problem is equivalent to a rectangular knapsack problem, and propose an MILP model and tailored genetic algorithms.
While many studies focus on static energy limits, 
\cite{homayouni2025optimizing} extend this by exploring speed-adjustable machines under peak power constraints in a job shop environment, minimizing both makespan and energy consumption. Similarly, \cite{hudaifah2025pareto} utilize a speed scaling strategy for just-in-time production systems to balance energy conservation with earliness and tardiness penalties, highlighting the trade-off between machine efficiency and schedule performance.

The previous discussion highlights the classical trade-off arising among execution times, resource (energy) consumptions and costs, among others.
Formulating bi-objective problems combining classical and energy-related objective functions is a standard way of approaching this aspect. However, since most mono-objective scheduling problems are $\mathcal{NP}$-hard \citep{pinedo2022scheduling}, their bi-objective counterparts are also $\mathcal{NP}$-hard and their Pareto frontiers are usually approximated through heuristic approaches. 
The majority of the recent literature in this field considers identical parallel machines. 
\cite{wang2018bi} address the BPMSTP using an augmented $\epsilon$-constrained method for small instances and an NSGA-II algorithm for medium and large scale instances.
\cite{ANGHINOLFI2021} propose an enhanced constructive heuristic combined with a novel Local Search (LS).
In a follow-up paper, \cite{GAGGERO2023} solve the BPMSTP by means of an exact solution algorithm based on the $\epsilon$-constrained method. One of their key contributions is the compact MILP formulation that exploits the structure of the solution space, reducing the computational burden of each single-objective subproblem and outperforming the results by \cite{ANGHINOLFI2021}.
The results by \cite{ANGHINOLFI2021} are also improved by \cite{JARBOUI2024}, who use a novel Variable Neighborhood Descent (VND) algorithm combined with dynamic programming to find the optimal energy cost on each machine, followed by an oscillating $\epsilon$-constraint method to optimize both TEC and makespan. We note that, to the best of our knowledge, no explicit comparison between \cite{JARBOUI2024} and \cite{GAGGERO2023} is available. 
Recently, \cite{LI2024} present a polynomial time approximation algorithm for the parallel machine scheduling problem with machine-dependent energy costs, where the goal is to minimize makespan subject to an upper bound on the TEC.

\cite{LI2022} study a bi-objective scheduling problem on uniform parallel batch processing machines minimizing the maximum lateness and the total pollution cost, and solve it with a discrete evolutionary algorithm. 
Regarding energy-efficient scheduling problems on non-identical parallel machines, \cite{LI2022} study a bi-objective scheduling problem on uniform parallel batch processing machines, minimizing the maximum lateness and the total pollution cost, and solve it with a discrete evolutionary algorithm.
Energy constraints have also been addressed for scheduling problems with unrelated machines (see, e.g., \citealt{fanjul2017models} for an overview of different variants and MILP models). 
\cite{CHOU2020} propose an original probabilistic energy usage model for unrelated parallel machines, introducing stochastic energy consumption functions for job/machine pairs from energy usage data over time. They propose a heuristic multi-objective algorithm to minimize energy costs and makespan under probabilistic power consumption and due date constraints.
\cite{cota2021} propose a matheuristic algorithm for a bi-objective variant considering sequence-dependent setup times, and minimizing makespan and total electricity consumption.
\cite{heydar2022approximate} propose an approximate dynamic programming algorithm for the unrelated parallel machines problem with incoming jobs at each time interval, where the objective is to minimize the expected value of a linear combination of the makespan and total energy cost.
\cite{WangChe2023} tackle a scheduling problem with TEC minimization with a bounded makespan using a Variable Neighborhood Search heuristic.
Some energy-aware problems also exist in the flow shop scheduling literature for energy saving in manufacturing (see, e.g., \citealt{duan2024}).

\subsection{Resource-constrained parallel machines}
The \problemname\ also belongs to the broad category of Resource-Constrained Parallel Machine Scheduling Problems (RCPMSP), where the goal is to schedule jobs across multiple machines while managing a set of resources required for processing each job.
The RCPMSP has been a subject of extensive research for decades \citep{BLAZEWICZ1983, EDIS2013}. 
Recent studies have explored various types of resources, including human workforce \citep{CASELLI2024} and automated guided vehicles \citep{Reddy2022}.
Additionally, the literature addresses both renewable  \citep{fanjul2017models, zhang2023improved} and non-renewable resources \citep{chen2024simplified}.
Resource-constrained scheduling also extends to the management of computational and network infrastructure. \cite{mendez2017} address energy-aware scheduling in homogeneous multicore systems, maximizing rewards from optional task execution through frequency scaling and migration policies. In the context of virtual data centers, \cite{belter2021scheduling} model bandwidth and physical links as critical resources, applying an on/off strategy to network equipment to minimize energy consumption during data transfers within fixed time windows.

\subsection{Time-dependent scheduling} 
Variable costs and task durations have also been studied in various optimization problems. 
Within the Vehicle Routing Problem (VRP) literature, variants that incorporate variable costs or travel times, usually dependent on the moment a given trip between two clients begins, are known as Time-Dependent VRPs (TDVRPs). Due to space limitations, we refer the reader to \cite{AdamoGGG24} for an updated survey on TDVRPs.
In the scheduling domain, significant attention has been given to problems with time-dependent processing times (e.g., deteriorating processing times).
These are categorized as \textit{time-dependent scheduling problems} \citep{gawiejnowicz2020, GAWIEJNOWICZ2025}. 
Recently, \cite{POTTEL2022} have addressed a single machine scheduling problem considering a predefined order of execution, a generic resource consumption function, and time-dependent activity durations.
This is referred to as the Time-Dependent Activity Scheduling Problem with resource constraints.
In some works, the processing time of jobs is a convex function of the amount of a non-renewable resource allocated
to the job, and the total amount of resource consumption is bounded (see, e.g., \citealt{CARRASCO2018}).
The \problemname\ addressed in this paper differs from the class of scheduling problems with time-dependent processing times, as those problems involve variable processing times but constant (fixed) consumption at each instant throughout their execution. In contrast, the PMSTVP focuses on scenarios where processing times are fixed but energy costs and resource consumption vary over time to represent realistic scenarios of scheduling under energy constraints. 

Different scheduling frameworks that do not involve machine assignment still address resource consumption. An important class of the problem is resource-constrained project scheduling \citep{ARTIGUES2025} with time-varying demands, where both activity requirements and resource availabilities may vary over time \citep{hartmann2013}. A further generalization is given by the continuous energy-constrained scheduling problem, a cumulative scheduling problem in which resource usage is modeled as a continuous, piecewise-constant function over time, constrained by concave processing functions and bounded energy requirements \citep{nattaf2017}. Similarly, time-dependent resource allocation problems extend the classical model of task assignment over time by allowing resource usage to vary dynamically during task execution \citep{PARK2024}.

From a methodological standpoint, the standard approach within the scheduling literature is to discretize the planning horizon and assume the resource availability and consumption (in our case, the energy) to be known and constant. 
The Time Indexed (or time-discrete) Formulations (TIFs, see, e.g., \citealt{KRAMER2021, SILVA2023}), where the decision variables are explicitly indexed by the (discrete) time horizon, stand as a convenient modeling framework for this type of problem. TIFs are flexible tools for modeling realistic variable resource consumption functions, although they can result in large models as the number of time periods increases, limiting the size of solvable instances.

Most of the papers discussed throughout this section, considering TOU tariffs, assume that the energy consumption of a job remains unchanged (i.e., constant and fixed) during its execution. The related models may represent a simplification of realistic situations.
In our case, considering variable energy consumption functions in addition to variable buying costs and selling prices is more realistic but introduces relevant additional challenges. 
From this review, we note that the \problemname\ has not been formalized within the energy-efficient, resource-constrained, or time-dependent scheduling literature, and that incorporating variable consumption within parallel machine scheduling motivates the design of new solution approaches to address the additional challenges.

\section{Problem definition}
\label{sec:problem}

Consider a set $J = \{1, \dots, n\}$ of jobs to be scheduled on a set $I = \{1, \dots, m\}$ of identical machines over a discrete, finite time horizon $T = \{0, \dots, \tmax\}$ during which the execution of all jobs must be completed. We refer to $t \in T$ as a time slot. Each job $j \in J$ has a processing time of $p_j$ time slots, during which energy is required for its execution. The energy consumption is machine-dependent and time-dependent, meaning that the energy consumed varies over time during the execution of the job. 
In accordance with the EES literature, we adopt the \textit{identical} parallel machine framework, in which job processing times are independent of the machine assignment (machines are identical in processing capability), while each machine has a different energy consumption rate (see, e.g., \citealt{CATANZARO2023, GAGGERO2023}).
The energy consumption function for job $j \in J$, when assigned to machine $i \in I$, is modeled as a step function (called function hereafter) $u_{ji\tau} \ge 0$ indicating the energy required at time $1 \leq \tau \leq p_j$ relative to the start of its execution.

Energy consumption is influenced by different characteristics of the infrastructure. 
We denote by $E$ the maximum amount of energy that the system can process to schedule jobs in every time slot $t \in T$. Energy can be obtained either from incoming clean, renewable sources (e.g., photovoltaic panels) or from the grid at a given cost.
Let $e_t \ge 0$ be the amount of energy obtained from the photovoltaic panels in time slot $t$. This energy varies depending on daylight and weather conditions and cannot be stored to be used in future time slots. 
Similar factors may influence the time-dependent availability of other renewable energy sources (e.g., wind or hydro), suggesting that the proposed model can be adapted to address analogous scheduling problems in different industrial settings.
At every time slot, the company may decide to buy (if needed) or sell energy (in case of excess). 
For each $t \in T$, let scalars $c_t, d_t \ge 0$ denote the cost of buying and the price of selling one unit of energy in time slot $t$, respectively. 
In addition, $c_t > d_t \; \forall t \in T$ holds and avoids the possibility of arbitrage. 
We assume that this information is deterministic and known in advance.
Modeling time-dependent buying and selling prices further supports the applicability of this model to a wide range of problems with dynamic or market-driven energy pricing systems.
The \problemname\ involves scheduling each job on exactly one machine, ensuring that no more than one job is executed on each machine at any given time slot. Additionally, the total energy consumed per time slot must not exceed $E$. The objective is to minimize the TEC, i.e., the total net energy cost of the schedule, defined as the difference between the total cost of energy acquisition and the revenue generated from selling energy.

The assumption that the energy consumption function $u_{ji\tau}$ is a piece-wise constant (step) function is directly motivated by the battery and electronic component testing operations described in Section~\ref{sec:intro}. These processes, which include repeated charge and discharge cycles, naturally result in highly irregular energy requirements and sudden power variations. By adopting a time-indexed representation of the energy consumption, we discretize the execution into time slots, allowing the model to capture these transitions while maintaining a mathematically tractable linear structure for the problem. Clearly, the quality of the approximation of the energy consumption depends on the granularity of the discretization. Finally, both TOU pricing and the energy generated from solar panels are handled in the same fashion, using step functions to represent their variability over the planning horizon.

Observe that the \problemname\ is $\mathcal{NP}$-hard because it generalizes the classical $P$$| |C_{\max}$ (\citealt{pinedo2022scheduling}). 
Moreover, finding a solution that satisfies the maximum energy constraint is as difficult as the subset sum problem (\citealt{CacchianiILM22}), and therefore even finding a feasible solution to the \problemname\ is $\mathcal{NP}$-complete.

%
We present now three examples to illustrate the impact of explicitly accounting for variable consumption functions in the \problemname, both with respect to the feasibility of a schedule and to its quality. 
In all examples, we consider $T = \{0,1,\dots,4\}$, assuming identical machines and $E = 4$. Each example is composed of two instances: (i) one with a variable consumption function for each job representing the real consumption, and (ii) another one with exactly the same parameters but with fixed consumption functions (i.e., $u_{ji\tau} = C$ for all $t, i$ for a given job $j$) with a total energy consumption that remains the same as in the realistic instance for all the jobs. In this way, the latter instance can be interpreted as a simplification of the former one.

\begin{example}
Figure~\ref{fig:1} shows an example with two jobs such that the \problemname\ instance with variable consumption functions is infeasible (Figures~\ref{fig:1b} and \ref{fig:1d}), but this becomes feasible if the jobs consider a fixed, average consumption throughout their execution (Figures~\ref{fig:1a} and \ref{fig:1c}). 
\rev{
The instances consider two jobs  with $p_{1} = p_{2} = 4$. Fig.~\ref{fig:1a} indicates the fixed consumption functions for each job $u_{1} = u_{2} = [2,2,2,2]$, while Fig.~\ref{fig:1b} depicts the variable consumption functions $u_{1} = [2,4,1,1]$ and $u_{2} = [1,2,4,1]$ for the realistic case. The schedule in Fig.~\ref{fig:1a} is feasible (Fig.~\ref{fig:1c}) while no feasible schedule for the variable consumption instance (Fig.~\ref{fig:1b}) exists because the energy budget $E = 4$ cannot be satisfied (Fig.~\ref{fig:1d}).}
Note that the instance with variable consumption is infeasible since any schedule would exceed the energy budget $E = 4$.

\begin{figure}[ht]
    \centering
    \begin{subfigure}[t]{0.45\textwidth}
    \centering
    \begin{tikzpicture}[scale=0.7]
        \draw (0,0) rectangle (5,2); 
        
        \foreach \x in {1,2,3,4} {
            \draw[gray, very thin] (\x,0) -- (\x,2);
        }
        \draw (0,1) -- (5,1);
        \node[left] at (-0.05,0.5) {\scriptsize \(m_1\)};
        \node[left] at (-0.05,1.5) {\scriptsize \(m_2\)};
        
        \fill[fill=gray!20] (0,0) rectangle (4,0.4); 

        \fill[fill=gray!50] (1,1) rectangle (5,1.4); 

        \node at (2,0.2) {\tiny \(1\)}; 
        \node at (3,1.2) {\tiny \(2\)}; 

        \foreach \x in {0,...,4} {
            \node[below, scale=0.7] at (\x+0.5,0) {\footnotesize \x};
        }
    \end{tikzpicture}
    \caption{Schedule with simplified jobs.}
    \label{fig:1a}
\end{subfigure}
    \hspace{-1cm}
    \begin{subfigure}[t]{0.45\textwidth}
        \centering
        \begin{tikzpicture}[scale=0.7]
            \draw (0,0) rectangle (5,2); 
            
            \foreach \x in {1,2,3,4} {
                \draw[gray, very thin] (\x,0) -- (\x,2);
            }
            \draw (0,1) -- (5,1);
            \node[left] at (-0.05,0.5) {\footnotesize \(m_1\)};
            \node[left] at (-0.05,1.5) {\footnotesize \(m_2\)};
            
            \fill[fill=gray!20] (0,0) rectangle (1,0.4); 
            \fill[fill=gray!20] (1,0) rectangle (2,0.8); 
            \fill[fill=gray!20] (2,0) rectangle (3,0.2); 
            \fill[fill=gray!20] (3,0) rectangle (4,0.2); 
    
            \fill[fill=gray!50] (1,1) rectangle (2,1.2); 
            \fill[fill=gray!50] (2,1) rectangle (3,1.4); 
            \fill[fill=gray!50] (3,1) rectangle (4,1.8); 
            \fill[fill=gray!50] (4,1) rectangle (5,1.2); 
    
            \node at (3.5, 1.3) {\tiny \(2\)}; 
            \node at (1.5, 0.3) {\tiny \(1\)}; 
    
            \foreach \x in {0,...,4} {
                \node[below, scale=0.7] at (\x+0.5,0) {\footnotesize \x};
            }
        \end{tikzpicture}
        \caption{Schedule with realistic jobs.}
        \label{fig:1b}
    \end{subfigure}
    %
    \begin{subfigure}[t]{0.45\textwidth}
            \centering
            \begin{tikzpicture}[scale=0.7]
                \draw (0,0) rectangle (5,2); 
                
                \foreach \x in {1,2,3,4} {
                    \draw[gray, very thin] (\x,0) -- (\x,2);
                }
                
                \draw[dashed] (0,0.8) -- (5,0.8);
                \node[left] at (-0.2, 0.8) {\footnotesize \(E\)}; 
        
        \draw[thick, black] 
            (0,0.4) -- (1,0.4) 
            -- (1,0.4) -- (1,0.8) 
            -- (2,0.8) -- (4,0.8) 
            -- (4,0.8) -- (4,0.4) 
            -- (4,0.4) -- (5,0.4); 

            \foreach \x in {0,...,4} {
                \node[below, scale=0.7] at (\x+0.5,0) {\footnotesize \x};
            }
            \end{tikzpicture}
            \caption{No energy limit violation.}
            \label{fig:1c}
        \end{subfigure}
    \hspace{-1cm}
    \begin{subfigure}[t]{0.45\textwidth}
        \centering
        \begin{tikzpicture}[scale=0.7]
            \draw (0,0) rectangle (5,2); 
            
            \foreach \x in {1,2,3,4} {
                \draw[gray, very thin] (\x,0) -- (\x,2);
            }
            
            \draw[dashed] (0,0.8) -- (5,0.8);
            \node[left] at (-0.2,1) {\footnotesize \(E\)}; 
    
    \draw[thick, black] 
        (0,0.4) -- (1,0.4) 
        -- (1,1.0) -- (2,1.0) 
        -- (2,0.6) -- (3,0.6) 
        -- (3,1.0) -- (4,1.0) 
        -- (4,0.2) -- (5,0.2); 

            \foreach \x in {0,...,4} {
                \node[below, scale=0.7] at (\x+0.5,0) {\footnotesize \x};
            }
        \end{tikzpicture}
        \caption{Energy limit violation.}
        \label{fig:1d}
    \end{subfigure}
    %
    \caption{\rev{Simplifying variable consumption to fixed consumptions may bring feasibility to an infeasible \problemname\ instance.}}
    \label{fig:1}
\end{figure}
\end{example}

\begin{example}
Conversely, Figure~\ref{fig:2} shows the opposite case. Consider again an example with two jobs but in this case the instance with variable consumption functions is feasible (Figures~\ref{fig:2b} and \ref{fig:2d}) while the corresponding instance with average fixed consumption functions is infeasible (Figures~\ref{fig:2a} and \ref{fig:2c}).
\rev{The instances contain two jobs with $p_{1} = p_{2} = 3$. Figure~\ref{fig:2a} shows the fixed functions $u_{1} = u_{2} = [3,3,3]$, while Figure~\ref{fig:2b} the variable functions $u_{1} = [4,4,1]$ and $u_{2} = [1,4,4]$ for the realistic case. The schedule presented in Figure~\ref{fig:2d} is feasible while any schedule for the fixed consumption instance is infeasible since it exceeds the energy budget $E = 4$ (Figure~\ref{fig:2c}).}
\begin{figure}[ht]
    \centering
    \begin{subfigure}[t]{0.45\textwidth}
    \centering
    \begin{tikzpicture}[scale=0.7]
        \draw (0,0) rectangle (5,2); 
        
        \foreach \x in {1,2,3,4} {
            \draw[gray, very thin] (\x,0) -- (\x,2);
        }
        \draw (0,1) -- (5,1);
        \node[left] at (-0.05,0.5) {\scriptsize \(m_1\)};
        \node[left] at (-0.05,1.5) {\scriptsize \(m_2\)};
        
        \fill[fill=gray!20] (2,0) rectangle (5,0.6); 

        \fill[fill=gray!50] (0,1) rectangle (3,1.6); 

        \node at (3.5,0.3) {\tiny \(1\)}; 
        \node at (1.5,1.3) {\tiny \(2\)}; 

        \foreach \x in {0,...,4} {
            \node[below, scale=0.7] at (\x+0.5,0) {\footnotesize \x};
        }
    \end{tikzpicture}
    \caption{Infeasible schedule for simplified jobs.}
    \label{fig:2a}
\end{subfigure}
    \hspace{-1cm}
    \begin{subfigure}[t]{0.45\textwidth}
        \centering
        \begin{tikzpicture}[scale=0.7]
            \draw (0,0) rectangle (5,2); 
            
            \foreach \x in {1,2,3,4} {
                \draw[gray, very thin] (\x,0) -- (\x,2);
            }
            \draw (0,1) -- (5,1);
            \node[left] at (-0.05,0.5) {\scriptsize \(m_1\)};
            \node[left] at (-0.05,1.5) {\scriptsize \(m_2\)};
            
            \fill[fill=gray!20] (2,0) rectangle (3,0.2); 
            \fill[fill=gray!20] (3,0) rectangle (4,0.8); 
            \fill[fill=gray!20] (4,0) rectangle (5,0.8); 
    
            \fill[fill=gray!50] (0,1) rectangle (1,1.8); 
            \fill[fill=gray!50] (1,1) rectangle (2,1.8); 
            \fill[fill=gray!50] (2,1) rectangle (3,1.2); 

            \node at (1.0, 1.4) {\tiny \(2\)}; 
            \node at (4.0, 0.4) {\tiny \(1\)}; 
    
            \foreach \x in {0,...,4} {
                \node[below, scale=0.7] at (\x+0.5,0) {\footnotesize \x};
            }
        \end{tikzpicture}
        \caption{Feasible schedule for realistic jobs.}
        \label{fig:2b}
    \end{subfigure}
    %
    \begin{subfigure}[t]{0.45\textwidth}
            \centering
            \begin{tikzpicture}[scale=0.7]
                \draw (0,0) rectangle (5,2); 
                
                \foreach \x in {1,2,3,4} {
                    \draw[gray, very thin] (\x,0) -- (\x,2);
                }
                
                \draw[dashed] (0,0.8) -- (5,0.8);
                \node[left] at (-0.2,0.8) {\footnotesize \(E\)}; 
        
        \draw[thick, black] 
            (0,0.6) -- (1,0.6) 
            -- (1,0.6) -- (2,0.6) 
            -- (2,1.2) -- (3,1.2) 
            -- (3,0.6) -- (4,0.6) 
            -- (4,0.6) -- (5,0.6); 

            \foreach \x in {0,...,4} {
                \node[below, scale=0.7] at (\x+0.5,0) {\footnotesize \x};
            }
            \end{tikzpicture}
            \caption{Energy limit violation.}
            \label{fig:2c}
        \end{subfigure}
    \hspace{-1cm}
    \begin{subfigure}[t]{0.45\textwidth}
        \centering
        \begin{tikzpicture}[scale=0.7]
            \draw (0,0) rectangle (5,2); 
            
            \foreach \x in {1,2,3,4} {
                \draw[gray, very thin] (\x,0) -- (\x,2);
            }
            
            \draw[dashed] (0,0.8) -- (5,0.8);
            \node[left] at (-0.2,0.8) {\footnotesize \(E\)}; 
    
    \draw[thick, black] 
        (0,0.8) -- (1,0.8) 
        -- (1,0.8) -- (2,0.8) 
        -- (2,0.4) -- (3,0.4) 
        -- (3,0.8) -- (4,0.8) 
        -- (4,0.8) -- (5,0.8); 

            \foreach \x in {0,...,4} {
                \node[below, scale=0.7] at (\x+0.5,0) {\footnotesize \x};
            }
        \end{tikzpicture}
        \caption{No energy limit violation.}
        \label{fig:2d}
    \end{subfigure}
    %
    \caption{\rev{Simplifying variable consumption to fixed consumptions may result in an infeasible \problemname\ instance.}}
    \label{fig:2}
\end{figure}
\end{example}

\begin{example}
Finally, we focus on the impact of the variable consumption functions in the objective function, in particular with respect to the TOU costs. Figure~\ref{fig:3} shows an example with only one job having $p_{1} = 3$, and consider the costs of buying energy $(c_t) = [0.10, 0.10, 0.01, 0.10, 0.10]$ and the prices of selling energy $d_t = 0$ for $t \in T$.
\begin{figure}[ht]
    \centering
    \begin{subfigure}[t]{0.45\textwidth}
    \centering
        \begin{tikzpicture}[scale=0.7]
            \draw (0,0) rectangle (5,2); 
            
            \foreach \x in {1,2,3,4} {
                \draw[gray, very thin] (\x,0) -- (\x,2);
            }
            \node[left] at (-0.05,1) {\scriptsize \(m_1\)};

            \node[left] at (-0.05,2.5) {\scriptsize \(c_t\)};
            \node at (0.5,2.5) {\tiny \(\textbf{0.10}\)};
            \node at (1.5,2.5) {\tiny \(\textbf{0.10}\)};
            \node at (2.5,2.5) {\tiny \(\textbf{0.01}\)};
            \node at (3.5,2.5) {\tiny \(\textbf{0.10}\)};
            \node at (4.5,2.5) {\tiny \(\textbf{0.10}\)};
            \node[right] at (5.05,2.5) {\tiny(\euro/kWh)};
            
            \fill[fill=gray!20] (1,0) rectangle (2,0.8); 
            \fill[fill=gray!20] (2,0) rectangle (3,0.8); 
            \fill[fill=gray!20] (3,0) rectangle (4,0.8); 

            \node at (2.5, 0.4) {\tiny \(1\)}; 
    
            \foreach \x in {0,...,4} {
                \node[below, scale=0.7] at (\x+0.5,0) {\footnotesize \x};
            }
        \end{tikzpicture}
    \caption{Simplified job with $z = 0.42$.}
    \label{fig:3a}
\end{subfigure}
    %
    \begin{subfigure}[t]{0.45\textwidth}
        \centering
        \begin{tikzpicture}[scale=0.7]
            \draw (0,0) rectangle (5,2); 
            
            \foreach \x in {1,2,3,4} {
                \draw[gray, very thin] (\x,0) -- (\x,2);
            }
            \node[left] at (-0.05,1) {\scriptsize \(m_1\)};

            \node[left] at (-0.05,2.5) {\scriptsize \(c_t\)};
            \node at (0.5,2.5) {\tiny \(\textbf{0.10}\)};
            \node at (1.5,2.5) {\tiny \(\textbf{0.10}\)};
            \node at (2.5,2.5) {\tiny \(\textbf{0.01}\)};
            \node at (3.5,2.5) {\tiny \(\textbf{0.10}\)};
            \node at (4.5,2.5) {\tiny \(\textbf{0.10}\)};
            \node[right] at (5.05,2.5) {\tiny(\euro/kWh)};
            
            \fill[fill=gray!20] (1,0) rectangle (2,0.4); 
            \fill[fill=gray!20] (2,0) rectangle (3,1.6); 
            \fill[fill=gray!20] (3,0) rectangle (4,0.4); 

            \node at (2.5, 0.4) {\tiny \(1\)}; 
    
            \foreach \x in {0,...,4} {
                \node[below, scale=0.7] at (\x+0.5,0) {\footnotesize \x};
            }
        \end{tikzpicture}
        \caption{Realistic job with $z = 0.24$.}
        \label{fig:3b}
    \end{subfigure}
    \caption{Simplifying realistic consumption may impact the total energy cost. 
    } 
    \label{fig:3}
\end{figure}
For a realistic consumption function $u_{1} = [1,4,1]$, the optimal solution (Figure~\ref{fig:3b}) has a cost of 0.24\euro. Any feasible solution for an average instance with $u_{1} = [2,2,2]$ (Figure~\ref{fig:3a}) has the cost of 0.42\euro. However, if each of these solutions is evaluated considering realistic energy consumption, only one of them matches the optimal objective value of 0.24\euro\ while the other two have an objective value of 0.51\euro. 
In sum, a solution for a realistic context obtained by assuming constant average energy consumption may lead to suboptimal solutions with significant room for improvement.
\end{example}

\section{Mixed Integer Linear Programming model}
\label{sec:model}

This section presents a TIF for the \problemname. Let $X_{jit} \in \{0,1\}$ be a binary variable taking the value one if and only if job $j \in J$ starts its execution on machine $i \in I$ in time slot $t \in T$. 
For any time slot $t \in T$, \rev{variables  $W_{jt}$ indicate the amount of energy consumed by job $j \in J$, and variables $U_t$ and $V_t$ indicate the amount of energy bought and sold, respectively}.
 
Let $T_{jt} = \{\max\{0,t-p_j+1\},...,t\}$ be the set of feasible starting times for job $j$ such that its execution spans until $t$. In other words, if job $j$ starts its execution in a time slot $\tau \in T_{jt}$, then its execution covers time slot $t \in T$. The MILP model reads:
\allowdisplaybreaks
\begin{align}
    \label{eq1_s4} & \min \quad \sum_{t \in T} (c_{t} U_t - d_{t} V_t)   \\
    \label{eq2_s4} & \text{s.t.} \quad \sum_{i \in I} \sum_{t=0}^{\tmax-p_j+1} X_{jit} = 1 \quad  && j \in J  \\
    \label{eq3_s4} & \quad \quad \sum_{j \in J} \sum_{\tau \in T_{jt}} X_{ji\tau} \le 1 \quad && i \in I, t \in T  \\
    \label{eq4_s4} & \quad \quad W_{jt} = \sum_{i \in I} \sum_{\tau \in T_{jt}} u_{ji(t-\tau+1)} X_{ji\tau} \quad && j \in J, t \in T \\
    \label{eq5_s4} & \quad \quad \sum_{j \in J} W_{jt} = e_t + U_t - V_t \quad && t \in T \\
    \label{eq6_s4} & \quad \quad e_t + U_t - V_t \le E \quad && t \in T \\
    \label{eq8_s4} & \quad \quad X_{jit} \in \{0,1\} \quad && j \in J, i \in I, t \in T, t \leq \tmax-p_j+1\\
    \label{eq9_s4} & \quad \quad \rev{W_{jt} \ge 0} \quad && j \in J, t \in T\\
    \label{eq10_s4} & \quad \quad \rev{U_t, V_t \ge 0} \quad && t \in T
\end{align}

The objective function \eqref{eq1_s4} minimizes the TEC (i.e., the difference between the total cost of buying and selling energy). Constraints \eqref{eq2_s4} guarantee that each job is scheduled exactly once. Constraints \eqref{eq3_s4} prevent the execution of more than one job in the same time slot on the same machine. Constraints \eqref{eq4_s4} link variables $X_{jit}$ and $W_{jt}$, computing the energy consumption of job $j$ in time slot $t$ based on its variable consumption function and its starting time. 
Constraints \eqref{eq5_s4} match the total energy consumed by the schedule with the net amount of energy processed in each time slot, while ensuring that there is no wasted energy (energy is bought from the grid if needed or sold to the grid if available in excess).
Constraints \eqref{eq6_s4} impose that the net amount of energy does not exceed the upper limit $E$ allowed by the system for executing the tests or, more generally, to execute the schedule.
As $c_t > d_t \; \forall t \in T$, no arbitrage is possible and any optimal solution satisfies either $U_t > 0$ or $V_t > 0$, but not both, in any time slot.
Constraints \eqref{eq8_s4}-\eqref{eq10_s4} define the domain of the variables.

\section{A matheuristic algorithm}
\label{sec:heur}

In this section, we present the building blocks of a matheuristic algorithm \citep{boschetti2022,boschetti2023} developed for tackling the \problemname. The approach, outlined in Algorithm~\ref{alg:ils}, follows the Iterated Local Search (ILS) paradigm \citep{Lourenco2010}, which has been applied to several scheduling problems (see, e.g., \citealt{coelho22025}). 

For the sake of notation, throughout this section we let $\mathcal{I} = (J, I, T, u = (u_{jit}), c = (c_t), d = (d_t), e = (e_t), E)$ denote an instance of the \problemname. Given a solution $s$, we let $(j,i,t) \in J \times I \times T$ indicate that job $j$ is scheduled on machine $i$ starting at time slot $t$. Aligned with standard scheduling terminology, $(j,i,t)$ is the schedule represented in solution $s$. Let also $f(s)$ denote the value of the objective function \eqref{eq1_s4} for solution $s$, with $f(s) = \infty$ if $s$ is infeasible. \rev{For a reference, $s_{\textrm{best}}$ tracks the best known feasible solution found during the execution.} We use total time and a maximum number of iterations without improvement as termination criteria. 

As noted in Section~\ref{sec:problem}, the feasibility problem for the \problemname\ is $\mathcal{NP}$-complete. To address this challenge, we employ a combination of constructive heuristics focused on feasibility, aiming to obtain an initial solution, which we name \mbox{\textsc{ConstructiveStep}}. \rev{One of these heuristics applies an adaptation of the MILP model described in Section~\ref{sec:model}. The specific details of \mbox{\textsc{ConstructiveStep}} are described in more detaile in Section~\ref{sec:heur:constructive}}.

Afterward, \rev{and starting from a feasible solution,} an improvement phase \mbox{\textsc{ImprovementStep}} combines a classical LS algorithm composed of \rev{a relocate and a swap operators, described in Section~\ref{sec:heur:LS}, and an operator built upon} the MILP model of Section~\ref{sec:model} \rev{to explore a larger neighbourhood, described in Section~\ref{sec:heur:milp-op}.} 

\rev{To escape from local optima, we employ} a perturbation scheme \rev{following the destroy-and-repair paradigm, abstracted within the \textsc{PerturbationStep} procedure,} that takes as input a parameter $\alpha \in (0,1)$ indicating the percentage of machines to be selected as input for the perturbation. \rev{This parameter $\alpha$ is adjusted depending on the evolution of the optimization, in particular with respect to the number of consecutive iterations without improvements. This step is enconded in the \textsc{PerturbationUpdate} procedure. These two procedures are described in more detail in Section~\ref{sec:heur:perturb}. For a reference, $k$ records the number of consecutive perturbation iterations without improvement.}

The computational complexity of Algorithm~\ref{alg:ils} is dominated by the number of times the \textsc{ImprovementStep} is executed. Since \textsc{ImprovementStep} considers tackling an MILP formulation for a restricted subproblem, its worst case complexity is exponential. For practical purposes, the implementation imposes a time limit for each call to the MILP.

\begin{algorithm}[ht]
\caption{ILS algorithm for the \problemname}
\label{alg:ils}
\begin{algorithmic}[1]
\State \rev{\textbf{Input:} $\mathcal{I} = (J, I, T, u, c, d, e, E)$.}
\State \rev{\textbf{Output:} $(s, f(s))$.}
\State $s_{\textrm{best}} \gets $\textsc{ConstructiveStep}($\mathcal{I}$)
\State \rev{$k \gets 1$}
\State \rev{$\alpha \gets $\textsc{PerturbationUpdate}$(k)$}\label{alg:ils:pert-init}
\State \textbf{if} $s_{\textrm{best}} = \{\}$, \textbf{then} \Return $(s_{\textrm{best}}, \infty)$
\While{termination criterion is not met}
    \State $s_{\textrm{new}} \gets $\textsc{ImprovementStep}$(\mathcal{I},s)$\label{alg:ils:improvement}
    \State \rev{$k \gets k + 1$}
    \If{$f(s_{\textrm{new}}) < f(s_{\textrm{best}})$}
        \State $s_{\textrm{best}} \gets s_{\textrm{new}}$
        \State \rev{$k \gets 1$}
    \EndIf
    \State \rev{$\alpha \gets $\textsc{PerturbationUpdate}$(k)$}\label{alg:ils:pert-update}
    \State $s \gets $\textsc{PerturbationStep}$(\mathcal{I},s_{\textrm{best}}, \alpha)$\label{alg:ils:perturbation}
\EndWhile
\State\Return $(s_{\textrm{best}}, f(s_{\textrm{best}}))$
\end{algorithmic}
\end{algorithm}

\subsection{Insertion heuristic}
We begin by describing a general insertion procedure that is used as a subroutine in different components of the \matheur. Algorithm~\ref{alg:insertion} depicts the pseudocode in a general fashion that can be adapted to follow different criteria regarding feasibility and optimality by slightly adjusting some specific steps and input parameters. 
\begin{algorithm}[ht]
\caption{Insertion heuristic}
\label{alg:insertion}
\begin{algorithmic}[1]
\State \rev{\textbf{Input:} $\mathcal{I} = (J, I, T, u, c, d, e, E)$; $\theta$; (restricted) feasible solution $\hat{s}$.}
\State \rev{\textbf{Output:} A feasible solution $s$ if possible; $\{\emptyset\}$ otherwise.}
\State \rev{$s \gets \hat{s}$}.
\State \rev{$J_0 \gets$ set of jobs $j \in J$ that are unassigned in $\hat{s}$}.
\State \rev{$L \gets$} sequence containing all jobs $j \in J_0$, sorted according to criterion $\theta$. \label{alg:insertion:step:sorting}
\While{$L$ is non-empty}
    \State Let $j$ be the next job in $L$. Remove $j$ from $L$.
    \State \parbox[t]{\dimexpr\textwidth-\leftmargin-\labelsep-\labelwidth}{For each $i \in I, t \in T$ check if it is feasible to schedule $j$ on $i$ starting at $t$ in $s$, that is, $t + p_j \le \tmax$, each time slot $\tau \in \{t,\dots, t + p_j\}$ is free on $i$ and $\omega_{\tau} + u_{ji(\tau - t)} \le E$.\strut}\label{alg:insertion:step:feasibility}
    \State \textbf{if} no feasible assignment is possible, \textbf{then} \Return \rev{$\{\}$}.
    \State Let $\hat{i} \in I, \hat{t} \in T$ be the selected machine and starting time slot, respectively.\label{alg:insertion:step:selection}
    \State \rev{\parbox[t]{\dimexpr\textwidth-\leftmargin-\labelsep-\labelwidth}{$s \gets s \cup \{(j,\hat{i},\hat{t})\}$.\strut}}
    \State \rev{\parbox[t]{\dimexpr\textwidth-\leftmargin-\labelsep-\labelwidth}{Mark $\tau$ as occupied on $\hat{i}$ for $\tau \in \{\hat{t},\dots,\hat{t}+p_j\}$\strut.}}
    \State \rev{\parbox[t]{\dimexpr\textwidth-\leftmargin-\labelsep-\labelwidth}{$\omega_{\tau} = {\omega}_{\tau} + u_{j\hat{i}(\tau - \hat{t})}$ for $\tau \in \{\hat{t},\dots,\hat{t}+p_j\}$\strut.}}
\EndWhile
\State\Return $s$
\end{algorithmic}
\end{algorithm}
The heuristic schedules one job at a time from a set $J_0$ of unassigned jobs, denoted by $L$, iteratively, according to the order defined by a selection criterion denoted by $\theta$. 
In our implementation, $\theta$ represents one of the following three orderings: increasing job index $j$, decreasing value of $p_j$, and random. We follow a \textit{first-fit} policy to accelerate execution. A \textit{best-fit} can be easily implemented by evaluating the objective function for the insertion of $j$ in each feasible pair $i \in I, t \in T$. However, limited preliminary experiments showed that this alternative was too expensive regarding computation times for a constructive heuristic, and therefore is disregarded as part of the implementation.
Among other constraints, feasibility tests require validating that the maximum amount of energy $E$ is not surpassed at any $t \in T$. Thus, let $\omega_t \ge 0$ denote the total energy consumption at $t \in T$ by the restricted solution $s$.

Note that by setting $\hat{s} = \{\emptyset\}$ this insertion heuristic becomes a constructive algorithm that attempts to generate an initial (feasible) solution. 
Let $\pmax = \max_{j \in J} p_j$. Steps~\ref{alg:insertion:step:feasibility}-\ref{alg:insertion:step:selection} can be implemented in $O(m\cdot\pmax\cdot\tmax)$ using standard data structures, since energy consumption is dependent on the machine assignment of jobs, and hence the worst-case complexity of the insertion heuristic is $O(n\cdot m \cdot \pmax \cdot \tmax)$, assuming this dominates the sorting in Step~\ref{alg:insertion:step:sorting}. We remark, however, that this complexity could be reduced by considering more sophisticated data structures. For instance, maintaining a list of \textit{maximal free intervals} for each machine representing blocks of consecutive free instants could significantly accelerate feasibility tests regarding the fit of the job in a given period. \rev{Furthermore, redundant structures such as segment trees supporting minimum and maximum range queries on the residual energy $r_t = E - \omega_t, t \in T$, can accelerate some feasibility checks. More precisely, let $\mu_{ji} = \max_{1 \le \tau \le p_j} u_{ji\tau}$ be the maximum consumption of a job $j \in J$ over machine $i \in I$. An interval $[\tau_1,\tau_2]$ can be directly discarded if (i) $\tau_2 - \tau_1 + 1 < p_j$, i.e. job $j$ does not fit, or (ii) $\max_{\tau_1 \le \tau \le \tau_2} r_{\tau} < \mu_{ji}$. Analogously, every free subinterval of length at least $p_j$ is feasible if $\min_{\tau_1 \le \tau \le \tau_2} r_{\tau} \ge \mu_{ji}$. Otherwise, a full search is still required.}

Let $F$ be the number of free segments in a machine. \rev{This list can be checked in $O(F)$, while the segment tree is built only once in $O(\tmax)$, and queried in $O(\log \tmax)$. In general, $F \le n_i \ll \tmax$, where $n_i$ stands for the number of jobs assigned to machine $i$ in the solution.} For highly fragmented solutions, which may be the case when time-dependent consumption and prices are considered, \rev{$F \approx n_i$. Considering that the insertion heuristic is not the bottleneck of the matheuristic, and that the running time of the ILS is dominated by the MILP search (see Section~\ref{sec:heur:milp-op}), the implementation of segment trees is not expected to produce relevant improvements. We do implement, however, the maximal free interval list per machine, seeking a trade-off between accelerating the feasibility tests and a more complex implementation. Note that this does not change the worst-case number of comparisons required by the feasibility check itself, which can still approach $\tmax$ in adversarial instances, but it substantially reduces the average number of candidates examined in practice}. 

\subsection{Perturbation}
\label{sec:heur:perturb}

To escape from local optima, \rev{the method \textsc{PerturbationStep}} (Algorithm~\ref{alg:ils}, Step~\ref{alg:ils:perturbation}) introduces a perturbation to the incumbent solution following a classical destroy-and-repair paradigm. Given a feasible solution $\hat{s}$ and a perturbation percentage $\alpha \in (0,1)$, we randomly select $\lceil \alpha m \rceil$ machines and remove all jobs assigned to them from $\hat{s}$, obtaining as a result a \textit{restricted feasible solution}. Afterwards, the repair step consists of executing the insertion heuristic (Algorithm ~\ref{alg:insertion}), \rev{where $\theta$ being set to a random order}.

The \rev{method \textsc{PerturbationUpdate} updates the} perturbation percentage $\alpha$ (Algorithm~\ref{alg:ils}, steps \ref{alg:ils:pert-init} and \ref{alg:ils:pert-update}) depending on the outcome of the improvement step (Algorithm~\ref{alg:ils}, Step~\ref{alg:ils:improvement}), \rev{more specifically depending on the number of consecutive iterations $k$ without improvements to the best known solution. The procedure is as follows.} Consider a sequence of $K$ increasing, non-overlapping intervals in $[0,1]$, that is, $(\underline{\beta}_k, \overline{\beta}_k) \subseteq [0,1]$, $\overline{\beta}_k \le \underline{\beta}_{k+1}$, $k = 1,\dots,K-1$.
Larger perturbations are considered as $k$ increases, meaning that more machines are selected for being emptied.
Initially, we set $k = 1$ and $\alpha$ is randomly selected from the smallest perturbation range $(\underline{\beta}_1, \bar{\beta}_1)$. Whenever an improving solution is found, we reset $k = 1$, restarting the process with $(\underline{\beta}_1, \bar{\beta}_1)$. Conversely, if the improvement phase is unable to improve $s_{\textrm{best}}$, we increase $k$ and choose $\alpha$ randomly from the next interval $k+1$, applying a larger perturbation. 

For example, let us assume \(m=10\) machines and three perturbation intervals $(\underline{\beta}_k,\bar{\beta}_k)$ as \((0.05,0.10)\), \((0.10,0.20)\), and \((0.20,0.35)\), where $K = 3$. Initially, \(k=1\) and \(\alpha\) is drawn from the first interval; if \(\alpha=0.08\), then \(\lceil \alpha m \rceil=1\) machine is selected and all its jobs are removed. If no improvement is found, \(k\) increases and \(\alpha\) is drawn from the next interval, leading, for instance, to removing \(2\) or \(3\) machines. Whenever an improving solution is obtained, \(k\) is reset to \(1\), returning to the smallest perturbation range.

\subsection{Improvement Step}
\label{sec:heur:improvement}
Starting from a feasible solution, the improvement phase of the \matheur\ (Algorithm~\ref{alg:ils}, Step~\ref{alg:ils:improvement}) combines two classical LS operators with an adaptation of the MILP model (Section~\ref{sec:model}) to be used as an improvement operator.

\subsubsection{Local Search}
\label{sec:heur:LS}

We begin by discussing the adaptation of two classical LS operators, \textit{swap} and \textit{relocate}, to the \problemname. The presence of arbitrary variable consumption functions introduces an additional complexity compared to the fixed case (see, e.g., \citealt{GAGGERO2023}), since the feasibility check requires an explicit validation of the energy consumed during the execution of the jobs. We perform two types of moves and accept only moves that keep the solution feasible. Namely: 

\begin{itemize}
    \item Swap: For each pair of different jobs $j_1, j_2 \in J$, $j_1, \ne j_2$, exchange the positions between $j_1$ and $j_2$, maintaining in each case the original starting time of the other job in each machine.
    More formally, let $(j_1, i_1, t_1), (j_2, i_2, t_2) \in s$, remove these two assignments, and attempt to schedule $(j_2, i_1, t_1)$ and $(j_1, i_2, t_2)$. Each of these moves can be evaluated in $O(\pmax)$, and therefore the overall neighborhood is explored in $O(n^2\cdot \pmax)$.
    \item Relocate: For each job $j \in J$, check if $j$ can be scheduled on a different time slot or on a different machine, without modifying the schedule of the remaining jobs in $J \textbackslash \{j\}$. Each possible move can be tested in $O(\pmax)$, and therefore the complete neighborhood can be explored in $O(n \cdot m \cdot \tmax \cdot \pmax)$.
\end{itemize}
\noindent Following the classification introduced in \cite{hansen2017variable}, our method implements a \textit{Basic VND}, i.e., it explores neighborhoods in sequential order and considers a \textit{best-improvement} selection criterion for each of them.

\subsubsection{MILP search}
\label{sec:heur:milp-op}
We consider a solution $s$ and a subset of machines $\tilde{I} \subseteq I$, and formulate an auxiliary problem to reallocate all jobs assigned to machines $i \in \tilde{I}$ in $s$ while maintaining the rest of the solution fixed. We explore this neighborhood using an adaptation of the MILP model \eqref{eq1_s4} - \eqref{eq10_s4}, where the variables corresponding to jobs in $I \setminus \tilde{I}$ remain fixed to their value in solution $s$, and tackle this restricted ILP model using a commercial solver for a limited time.
This procedure is repeated several times within the MILP search, changing both the size and the composition of $\tilde{I}$, until an improving solution is found or a maximum number of iterations is reached.

By design, this procedure is more time-consuming than the LS (as it explores a larger neighborhood), and it is invoked only in case the former fails in finding an improving solution. Algorithm~\ref{alg:milp-op} depicts a pseudocode of the MILP search.
\begin{algorithm}
\caption{MILP search}
\label{alg:milp-op}
\begin{algorithmic}[1]
\State \rev{\textbf{Input:} $\mathcal{I} = (J, I, T, u, c, d, e, E)$; a feasible solution $s$.}
\State \rev{\textbf{Output:} A feasible solution $s$.}
\State \rev{$k \gets 1$.}
\State \rev{$\alpha \gets $\textsc{PerturbationUpdate}$^{*}(k)$.}\label{alg:milp-op:pert-init}
\While{$s$ not improved and maximum iterations not reached}
    \State Randomly select a subset $\tilde{I} \subseteq I$ with \rev{$\alpha$} percentage of the machines
    \State \parbox[t]{\dimexpr\textwidth-\leftmargin-\labelsep-\labelwidth}{Let $\tilde{J} = \{j \in J : (j,i,t) \in s, i \in \tilde{I}\}$.\strut}
    \State \rev{\parbox[t]{\dimexpr\textwidth-\leftmargin-\labelsep-\labelwidth}{$\tilde{s} \gets$ restricted solution by removing from $s$ the assignments $(j,i,t) \in s, j \in \tilde{J}$\strut}}
    \State \rev{\parbox[t]{\dimexpr\textwidth-\leftmargin-\labelsep-\labelwidth}{Formulate model \eqref{eq1_s4} - \eqref{eq10_s4}, fix the assignment of jobs $j \in J \setminus \tilde{J}$, i.e., set $X_{jit} = 1$ for $(j,i,t) \in \tilde{s}$.\strut}}\label{alg:milp-op:step:repair-model}
    \State \parbox[t]{\dimexpr\textwidth-\leftmargin-\labelsep-\labelwidth}{$s^* \gets$ solution obtained through an MILP solver with a time limit of $\gamma$ seconds..\strut}\label{alg:milp-op:step:repair}
    \If{$f(s{^*}) < f(s)$}
        \State $s \gets s^*$
    \EndIf
    \State \rev{$k \gets k + 1$.}
    \State \rev{$\alpha \gets$ \textsc{PerturbationUpdate}$^*(k)$.} \label{alg:milp:update-pert}
    \State Update time limit $\gamma$. \label{alg:milp:update-time}
\EndWhile
\State\Return $s$
\end{algorithmic}
\end{algorithm}
It also follows a destroy-and-repair scheme, similar to the perturbation described in Section~\ref{sec:heur:perturb}. For the repair step, a few observations are worth discussing. Step~\ref{alg:milp-op:step:repair} considers for each time slot the remaining energy budget at time $t \in T$ for the schedule, that is, $\hat{E}_t = E - \hat{\omega}_t$, and recall that $\hat{\omega}_t \ge 0$ denotes the total energy consumption at $t$ by the restricted solution $\hat{s}$. 
Furthermore, the restricted MILP can reschedule the jobs $j \in \tilde{J}$ that were removed from solution $s$ to any machine with an available and feasible time slot, rather than limiting them to the original subset of machines $\tilde{I}$.
From preliminary computational experiments, we observed that this general approach provides better solutions compared to its restricted variant in which we reallocate jobs $j \in \tilde{J}$ only to machines in $\tilde{I}$, even when the resulting MILP model may become larger and, consequently, more difficult to solve. 
Finally, the perturbation step in the MILP search follows a similar approach to that in ILS. \rev{The only difference relies in the update of the perturbation percentage, where \textsc{PerturbationUpdate}$^*$ (Step~\ref{alg:milp:update-pert}) updates parameter $\alpha$ only after five consecutive unsuccessful attempts to improve the current solution using the same value of $\alpha$}.
For larger problem instances, values of $\alpha$ close to 1 are limited to prevent the MILP from growing excessively large, ensuring it remains computationally manageable within the ILS framework.
In Step~\ref{alg:milp:update-time}, the time limit $\gamma$ for solving the adapted MILP is increased, allowing the model additional time to explore improved solutions.
In the next section, we introduce a further adaptation of the MILP search used in \matheur\ to build an initial feasible solution.

\subsection{Constructive Step}
\label{sec:heur:constructive}

Finding a feasible solution is one of the critical steps in our solution approach, both for the importance within the ILS approach and the complexity of the feasibility problem of the \problemname. For this reason, we attempt a sequence of steps to overcome this aspect, in many cases reusing or adapting some of the components described previously. The overall approach is as follows:

\begin{enumerate}
    \item We run the insertion heuristic (Algorithm \ref{alg:insertion}) with \rev{$\hat{s} = \{\}$} (all jobs to be inserted) multiple times, changing each time the sorting criteria $\theta$ to increasing job index, decreasing $p_j$, increasing $p_j$, and to a random order for 1000 iterations.
    \item If no feasible solution is returned, the last generated random sequence is returned to the insertion heuristic relaxing energy budget constraint.
    The adapted MILP search takes an infeasible solution and applies the same procedure as in Algorithm \ref{alg:milp-op} but taking the first feasible solution as stopping criterion.
    \item If Step 2 fails, we attempt to solve the complete model \eqref{eq1_s4} - \eqref{eq10_s4} until a feasible (not necessarily optimal) solution is found, if any.
\end{enumerate}  

\section{Computational experiments}
\label{sec:exp}
This section presents the outcome of extensive computational tests aimed at assessing the importance of incorporating variable consumption functions in the \problemname\ and evaluating the performance of the algorithms. 
Section~\ref{subsec:instances} describes the experimental setup and the generation of \problemname\ instances.
\rev{Section~\ref{subsec:sa-benchmark} briefly describes the SA metaheuristic benchmark developed to compare against our \matheur.}
Subsequently, Section~\ref{subsec:exp} reports a series of computational tests. 
First, in Section~\ref{subsec:res-comp}, we analyze the contribution of each component of \matheur\ proposed in Section~\ref{sec:heur}. 
Then, in Section~\ref{sec:overall-comparison},  we compare \matheur\ with the \milpmodel\ model \eqref{eq1_s4} - \eqref{eq10_s4} \rev{and SA} on small and medium instances. Finally\rev{,} in Section~\ref{subsec:res-larger}, we evaluate the performance of \matheur\ against SA on large instances. 
To conclude, in Section~\ref{subsec:qualitative-results} we analyze the solutions \rev{to highlight the challenges of our parallel machine scheduling problem and the impact of variable consumption on instance
feasibility.}

The experiments were executed on a single thread of a virtual machine Intel(R) Xeon(R) Gold with 2.30 GHz and 16 GB of RAM, running under Windows 10 Pro N. 
All algorithms were coded in \texttt{Python}, and \texttt{Gurobi 11} was used as the MILP solver. To accelerate the execution of the \matheur\ code, we resorted to the \texttt{Numba} library to precompile the most frequently used functions of the code and, therefore, improve the computation times. 

The methods were tested on sets of instances generated in order to assess the impact of variable energy consumption functions and described in Section~\ref{subsec:instances}.  Briefly, each instance with variable energy consumption is paired with a corresponding instance having the same parameters except for the energy consumption of the jobs, which is assumed to follow a fixed consumption for each job but maintaining the same total energy usage over its execution.
Results on related fixed-consumption instances are reported and commented in Appendix \ref{App:A}. To foster reproducibility, all instances and their solutions can be retrieved from \href{https://github.com/regor-unimore/Energy-Efficient-Scheduling-Variable-Consumption.git}{https://github.com/regor-unimore/Energy-Efficient-Scheduling-Variable-Consumption.git}, which also includes a graphical solution visualizer.
For \milpmodel, we set a time limit of 1200 seconds, while for \matheur\ and SA, this limit is reduced to 600 seconds to consider a more restricted context for their execution.

\subsection{Experimental setup}
\label{subsec:instances}
The instances for the \problemname\ are determined by:
\begin{itemize}
    \item the number $|J|$ of jobs, the number $|I|$ of machines, and the length $|T|$ of the planning horizon: we initially test 9 combinations $(|J|,|I|)$ $\in$ $\{(5,3)$, $(15,5)$, $(20,7)$, $(30,10)$, $(60,15)$, $(90,20)$, $(100,25)$, $(150,30)$, $(200,35)\}$ and 4 combinations $|T| \in \{48,72,96,120\}$ for small- and medium-size instances;
    \item an additional set of large-size instances determined by one size combination $(|J|,|I|, |T|)$ $=$ $(300,45,144)$;
    \item the instance \textit{saturation} $\nu$, that is the ratio between total processing time to be scheduled over all jobs and total available scheduling time from machines: we consider $\nu$ $\in \{0.7, 0.8, 0.9\}$;
    \item the processing time $p_j$ of a job $j$, that is randomly sampled from a normal distribution $N(\mu, \frac{\mu}{3})$, truncated to zero, where $\mu = \frac{|T|\cdot|I|}{|J|}\nu$ is the average processing time per job required for a given saturation factor $\nu$;
    \item the energy consumption: each machine $i \in I$ is randomly assigned to one of three levels of average energy consumption per hour equal to 30, 50 or 70 kWh to represent small, medium and large machines, respectively.
    We set the energy budget per time slot $E$ as the sum of energy levels of all the machines rescaled according to the time duration of the time slots obtaining values in kW per time slot.
\end{itemize}
Time-dependent information is generated as follows. The costs to buy energy are set according to the standard contracts in the Italian market. A TOU tariff policy with variable costs is derived by assigning different tiers depending on the moment of the day (Table~\ref{table:tou-tariffs}). We assume the same structure for the prices to sell energy but set the value to one third of the buying cost. Likewise, we assume an ideal production of a photovoltaic system by a step function (Figure~\ref{fig:solar-panels}). Recalling that $T$ represents 24 hours, values for $c_t, d_t$ and $e_t$ are defined by simply matching time slots $t \in T$ with its corresponding time duration.

\begin{minipage}{0.3\textwidth}
    \centering
    \captionof{table}{TOU-tariff scheme.}
    \scriptsize
    \begin{tabular}{lr}
        \toprule
        Time & Energy cost (\euro / kWh)\\
        \midrule
        00:00 - 07:00 & 0.12 \\
        07:00 - 08:00 & 0.15 \\
        08:00 - 19:00 & 0.18 \\
        19:00 - 23:00 & 0.15 \\
        23:00 - 23:59 & 0.12 \\
        \bottomrule
    \end{tabular}
   
    \label{table:tou-tariffs}
\end{minipage}
\hfill
\begin{minipage}{0.5\textwidth}
    \vspace{1.3cm}
    \centering
    \begin{tikzpicture}[scale=0.30, transform shape]
        \tikzstyle{axislabel}=[font=\Huge]
        
        \draw[thick, ->] (0,0) -- (24.5,0) node[anchor=north] {};
        \draw[thick, ->] (0,0) -- (0,8) node[left, xshift=-15pt] {\Huge kWh};
        \foreach \x in {0,4,8,12,16,20,24} {
            \draw (\x,0.1) -- (\x,-0.3) node[anchor=north] {\Huge \x:00};
        }
        \draw (-0.3,3.125) -- (-0.1,3.125) node[anchor=east] {\Huge 125};
        \draw (-0.3,6.25) -- (-0.1,6.25) node[anchor=east] {\Huge 250};
        \fill[gray!30]
         (0,0) -- (8,0) 
        -- (8,1.04) -- (9,1.04) 
        -- (9,2.08) -- (10,2.08) 
        -- (10,3.13) -- (11,3.13) 
        -- (11,4.17) -- (12,4.17) 
        -- (12,5.21) -- (13,5.21) 
        -- (13,6.25) -- (14,6.25) 
        -- (14,5.21) -- (15,5.21) 
        -- (15,4.17) -- (16,4.17) 
        -- (16,3.13) -- (17,3.13)
        -- (17,2.08) -- (18,2.08)
        -- (18,1.04) -- (19,1.04) 
        -- (19,0) -- (24,0) -- cycle ;
        
        \draw[thick] 
        (0,0) -- (8,0) 
        -- (8,1.04) -- (9,1.04) 
        -- (9,2.08) -- (10,2.08) 
        -- (10,3.13) -- (11,3.13) 
        -- (11,4.17) -- (12,4.17) 
        -- (12,5.21) -- (13,5.21) 
        -- (13,6.25) -- (14,6.25) 
        -- (14,5.21) -- (15,5.21) 
        -- (15,4.17) -- (16,4.17) 
        -- (16,3.13) -- (17,3.13)
        -- (17,2.08) -- (18,2.08)
        -- (18,1.04) -- (19,1.04) 
        -- (19,0) -- (24,0) ;
    \end{tikzpicture}
    
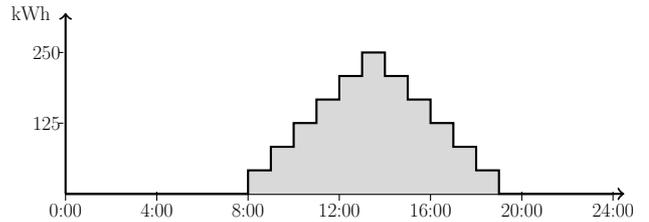
\captionof{figure}{Photovoltaic energy production.}
    \label{fig:solar-panels}
\end{minipage}%
\vspace{0.5cm}

\noindent From these definitions we obtain $9 \times 4 \times 3 = 108$ \textit{base configurations} for small-to-medium instances, and $45$ additional base configurations for the set of large instances. An instance for the \problemname\ is defined by a combination of a base configuration and a general consumption function. For each $j \in J$ and $0 \le \tau \le p_j$ we define a consumption function $v_{j\tau} \ge 0$ to be machine independent, and then incorporate the dependence on the machine by multiplying the machine level of consumption and obtain $u_{ji\tau}$ (in kW per time slot). 
We generate two types of instances, using two different energy consumption functions, as follows:

\begin{itemize}
    \item \textit{Fixed consumption} represents constant energy consumption, i.e., $v_{j\tau} = 1$ for $\tau \in  [1,p_j]$, obtaining a total area (energy consumption) of $p_j$ for each job $j$.
    \item \textit{Variable consumption} represents \textit{realistic} energy consumptions. For each job $j$ with processing time $p_j$, values $v_{j\tau}$ are generated randomly in a controlled fashion, guaranteeing that the total energy consumption is $p_j$. 
    For $\tau \in  [1,p_{j}-1]$, let $\bar{v}_{\tau} = (p_j - \sum\limits_{r = 0}^{\tau-1} v_{jr})/(p_j - \tau +1)$ be a measure of the average consumption for the subsequent time slots $\tau+1,\dots,p_j$, assuming that $v_{j0}=0$ and the overall consumption must be $p_j$ (as in the fixed consumption). Then, we iteratively set $v_{j\tau} \in \textrm{Unif}(\frac{1}{2}\bar{v}_{\tau},\frac{3}{2}\bar{v}_{\tau})$. For $\tau = p_j$, we set $v_{jp_j} = p_j - \sum\limits_{\tau = 1}^{p_j-1} v_{j\tau}$.
\end{itemize}

\noindent 
The above setup guarantees that the total energy consumption of a job $j \in J$ in a base configuration is the same across the two derived instances, i.e., $\sum\limits_{\tau = 1}^{p_j}u_{ji\tau}$ is exactly the same as given $i \in I$ for the corresponding fixed and variable consumption instances. By design, instances with a fixed consumption can be interpreted as a simplification of the variable ones in each base configuration. 
Moreover, schedules that fit within the planning horizon are always feasible for instances with fixed consumption, while this is not necessarily the case for the realistic consumption patterns. 

\newcommand{\ninst}{\#inst}
\newcommand{\nfeas}{\#feas}
\newcommand{\nopt}{\#opt}
\newcommand{\fG}{\%fG}
\newcommand{\comptime}{time}
\newcommand{\inctime}{inc.time}
\newcommand{\iG}{\%impG}
\newcommand{\avglb}{Avg. LB}
\newcommand{\avgub}{Avg. UB}
\newcommand{\fg}{\%gap}
\newcommand{\fallg}{\%fALLg}
\newcommand{\impg}{\%impg}

\subsection{\rev{The Simulated Annealing metaheuristic benchmark}}
\label{subsec:sa-benchmark}

\rev{As a benchmark, we compare our approach with a standard SA algorithm, commonly used for complex machine scheduling problems (see, e.g., \citealt{bektur2019mathematical,LiEtAl2024}). The SA neighborhood combines swap and relocate moves, randomly selected and applied at each iteration, and relies on the same constructive scheme as the \matheur\ considered in these experiments.

The SA cooling schedule is governed by three parameters: the initial temperature $t_0$, the cooling rate $\alpha$, and the number of iterations $\Delta_k$ performed at each temperature level before the temperature is decreased. These parameters were calibrated using the \texttt{irace} package \citep{Lopez2016} over a representative subset of the tested instances, with a budget of 200 experiments (\texttt{maxExperiments}) and a time limit of 200 seconds per run during the tuning phase. The following parameter ranges were considered:
\begin{itemize}
    \item $t_0 \in \{1, 5, 10, 25, 50, 75, 100, 250, 500\}$;
    \item $\alpha \in \{0.70, 0.80, 0.85, 0.90, 0.95, 0.99, 0.995\}$;
    \item $\Delta_k \in \{1, 5, 10, 50, 100, 250, 500\}$.
\end{itemize}
The best configuration returned by \texttt{irace} is $t_0 = 50$, $\alpha = 0.95$, and $\Delta_k = 50$; we remark that none of these values lies at the boundary of the corresponding range, indicating that the considered search space was adequately specified. The resulting SA configuration is used throughout the experiments in Sections~\ref{sec:overall-comparison} and~\ref{subsec:res-larger}.}

\subsection{Evaluation of algorithms}
\label{subsec:exp}
In this section, we delve into the experimental results focusing on the algorithmic performance.
\subsubsection{Analysis of \matheur\ components}
\label{subsec:res-comp}
First, we analyze in detail the performance of the \matheur\ presented in Section \ref{sec:heur} on small-to-medium-sized instances, focusing on the improvement step by isolating the contribution of each LS operator within \matheur, namely swap, relocate, MILP search, to the final performance. Table \ref{table:results-components} reports the results of one run of execution over all instances with variable consumption functions, first grouped according to the number $|T|$ of time slots, and then for each group we report the aggregated results according to the number of jobs and machines $(|J|$, $|I|)$. For each group, \ninst\ indicates the total number of instances within the group. Then, for each combination of \matheur\ components described in Section \ref{sec:heur:LS}, the table reports the number of optimal solutions found (\nopt), 
the average incumbent time (\inctime), and the average improvement gaps with respect to MILP (\%gapUB). 
The latter measures the average improvement gap with respect to \milpmodel, calculated as $100\times(z^* - z_{\textrm{ub}})/z^*$, with $z_{\textrm{ub}}$ and $z^*$ be the value of the best incumbent solution of, respectively, the \milpmodel\ model and the selected combination of \matheur\ components.
This gap is computed only on instances for which a feasible solution of MILP exists (indicated as \#feas in column ``MILP'' on the left).
Note that a negative value of \%gapUB indicates an improvement with respect to the \milpmodel\  solution.
The number of feasible solutions (\nfeas) for all the ILS variants is reported in a unique column Heur, as it only depends on the constructive step and remains unaffected when modifying the local search components of \matheur.

\begin{table}[ht]
\centering
\caption{Analysis of the impact of the LS operators in \matheur. Values of \inctime\ are expressed in seconds.}
\resizebox{\textwidth}{!}{%
\begin{tabular}{llc rr rrr rrr rrr rrr}
\toprule
\multicolumn{3}{c}{Instances} & \multicolumn{1}{c}{{MILP}} & \multicolumn{1}{c}{{Heur}} & \multicolumn{3}{c}{swap + relocate} &  \multicolumn{3}{c}{relocate + MILP search} &  \multicolumn{3}{c}{swap + MILP search} &  \multicolumn{3}{c}{swap + relocate + MILP search} \\
\cmidrule(lr){1-3} \cmidrule(lr){4-4}\cmidrule(lr){5-5} \cmidrule(lr){6-8} \cmidrule(lr){9-11} \cmidrule(lr){12-14} \cmidrule(lr){15-17} 
$|T|$ & $(|J|$, $|I|)$ & \ninst & {\#feas} & {\#feas} & \nopt  & \inctime & \%gapUB  & \nopt  & \inctime & \%gapUB  & \nopt  & \inctime & \%gapUB  & \nopt  & \inctime & \%gapUB  \\
\cmidrule(lr){1-3} \cmidrule(lr){4-4}\cmidrule(lr){5-5}  \cmidrule(lr){6-8} \cmidrule(lr){9-11} \cmidrule(lr){12-14} \cmidrule(lr){15-17} 
\multirow{3}*{48} 
& ([5, 20], [3, 7]) & 9 & {7} & {7} &	3	&	134.5	&	6.00	&	4		&	276.1	&	-0.08	&	4	&	287.7	&	0.01	& 4 & 220.3 & 0.09 \\
& ([30, 90], [10, 20])  & 9   & {9} & {9} &	0		&	48.7	&	2.93	&	0		&	527.9	&	0.27	&	0	&		425.1	&	0.22	& 0  & 490.0 & 0.24 \\
& ([100, 200], [25, 35])  & 9 & {9} & {9} &	0		&	49.7	&	1.25	&	0		&	588.3	&	1.98	&	0		&	587.5 &	0.74 & 0  & 481.2 & 0.63 \\
\cmidrule(lr){2-3}\cmidrule(lr){4-4}\cmidrule(lr){5-5}  \cmidrule(lr){6-8} \cmidrule(lr){9-11} \cmidrule(lr){12-14} \cmidrule(lr){15-17} 
&  & 27 & {25} & {25} &	3		&	77.7	&	3.18	&	4	&		464.1	&	0.79	&	4		&	433.4	&	0.35	& 4 &  397.2 & 0.34 \\
\cmidrule(lr){1-3} \cmidrule(lr){4-4}\cmidrule(lr){5-5}  \cmidrule(lr){6-8} \cmidrule(lr){9-11} \cmidrule(lr){12-14} \cmidrule(lr){15-17} 
\multirow{3}*{72}
& ([5, 20], [3, 7]) & 9& {8} & {7} &	1		&	137.4	&	4.42	&	3		&	299.5	&	0.16	&	3		&	229.4	&	0.16 & 3  & 302.2 &  0.11 \\
& ([30, 90], [10, 20])  & 9  & {9} & {9} &	0		&	7.9	&	4.55	&	0		&	531.6	& 0.73	&	0		&	545.3	&	0.44	& 0  & 560.2 & 0.36 \\
& ([100, 200], [25, 35])  & 9 & {7} & {9} &	0		&	68.4	&	-4.02	&	0		&	576.3	&	-2.58	&	0	&		582.8	&	-3.49	& 0  & 574.5 & -4.73  \\
\cmidrule(lr){2-3}\cmidrule(lr){4-4}\cmidrule(lr){5-5} \cmidrule(lr){6-8} \cmidrule(lr){9-11} \cmidrule(lr){12-14} \cmidrule(lr){15-17} 
&  & 27  & {24} & {25} &	1		&	71.2	&	2.01	&	3	&		469.1	&	-0.43	&	3		&	452.5	& -0.80	& 3 & 478.9 &  1.21 \\
\cmidrule(lr){1-3}\cmidrule(lr){4-4}\cmidrule(lr){5-5} \cmidrule(lr){6-8} \cmidrule(lr){9-11} \cmidrule(lr){12-14} \cmidrule(lr){15-17} 
\multirow{3}*{96} 
& ([5, 20], [3, 7]) & 9 & {6} & {6} &	1		&	247.8	&	5.97	&	3		&	322.9	&	0.17	&	3		&	332.0	&	-0.07	& 3 & 339.7 &  -0.02 \\
& ([30, 90], [10, 20])  & 9 & {8} & {9} &	0		&	35.3	&	4.13	&	0		&	533.6	&	0.81	&	0	&		554.7	&	0.28	& 0 &  550.1 &  0.57\\
& ([100, 200], [25, 35])  & 9 & {6} & {9} &	0		&	79.6	&	-16.03	&	0		&	587.2	&	-14.08	&	0	&		575.7	&	-14.87	& 0  & 553.3 & -16.70  \\
\cmidrule(lr){2-3} \cmidrule(lr){4-4}\cmidrule(lr){5-5}  \cmidrule(lr){6-8} \cmidrule(lr){9-11} \cmidrule(lr){12-14} \cmidrule(lr){15-17} 
&  & 27  & {20} & {24} &	1		&	120.9	&	-1.37	&	3	&		481.3	&	-3.85	&	3		&	487.5	&	-4.37	& 3   & 481.0 & -4.79 \\
\cmidrule(lr){1-3} \cmidrule(lr){4-4}\cmidrule(lr){5-5} \cmidrule(lr){6-8} \cmidrule(lr){9-11} \cmidrule(lr){12-14} \cmidrule(lr){15-17} 
\multirow{3}*{120} 
& ([5, 20], [3, 7]) & 9  & {5} & {5} & 2	&	331.6	&	5.74	&	3		&	422.8	&	0.45	&	3		&	422.0	&	0.55	& 3   & 464.2 & 0.59  \\
& ([30, 90], [10, 20])  & 9  & {5} & {7} &	0		&	163.9	& 4.36	&	0		&	581.1	&	1.61	&	0	&		568.5	&	0.60	& 0   & 556.5 & 1.11 \\
& ([100, 200], [25, 35])  & 9 & {3} & {7} &	0		&	226.1	&	-15.84	&	0	&	570.5	&	-12.84	&	0	&		583.4	&	-10.98	& 0  & 531.5 & -16.06  \\
\cmidrule(lr){2-3} \cmidrule(lr){4-4}\cmidrule(lr){5-5}  \cmidrule(lr){6-8} \cmidrule(lr){9-11} \cmidrule(lr){12-14} \cmidrule(lr){15-17} 
&  & 27  & {13} & {19} & 2 & 240.6 & 0.23 & 3 & 524.8 & -2.17 & 3 & 524.6 & -2.09 & 3 & 517.4 & -3.05 \\
\cmidrule(lr){1-3} \cmidrule(lr){4-4}\cmidrule(lr){5-5} \cmidrule(lr){6-8} \cmidrule(lr){9-11} \cmidrule(lr){12-14} \cmidrule(lr){15-17} 
Tot &  & 108 & {82} & {93} &	7	&	127.6	&	1.26	&	13	&		484.8	&	-1.17	&	13		&	474.5	& -1.53	& 13 &  468.6 &  -1.90 \\
\bottomrule
\end{tabular}%
}
\label{table:results-components}
\end{table}
The results confirm that the LS operators within \matheur\ contribute to the performance of the algorithm, especially the MILP search (see Section \ref{sec:heur:milp-op}), which is able to explore larger neighborhoods. Indeed, removing the MILP search (swap + relocate combination) reduces average incumbent times, although the quality of the solutions degrades significantly (this was noticed also when increasing by ten times the maximum number of iterations without improvements).
Excluding the MILP search increases the average improvement gap on all 108 instances by 3.19\% (from 1.26\% to -1.9\%).
These results validate the effectiveness of our \matheur\ algorithm, particularly the matheuristic design incorporating the MILP search.

\subsubsection{Comparison of \matheur\rev{, SA,} and \milpmodel}
\label{sec:overall-comparison}
After assessing the validity of our \matheur\ algorithm, we compare it with state-of-the-art alternative solution methods.
Table \ref{table:results-random-summary} presents the results for \milpmodel\rev{, SA, and \matheur}, grouped in a similar manner to Table \ref{table:results-components}. Here, for each method, \nfeas\ and \nopt\ indicate the number of instances for which the corresponding method is able to find a feasible and an optimal solution, respectively, within the time limit. 
For \milpmodel, let $z_{\textrm{ub}}$ and $z_{\textrm{lb}}$ denote the objective value of the incumbent solution and the best bound obtained, respectively. We denote by \%gapLB the average final gap, calculated as $100\times ($$z_{\textrm{ub}}$ - $z_{\textrm{lb}})/z_{\textrm{ub}}$ for each instance if $z_{\textrm{lb}}$ and $z_{\textrm{ub}}$ are known, and by \comptime\ the average computation time required by \texttt{Gurobi}. \rev{The average values of the lower and upper bounds obtained for each group of instances, \avglb\ and \avgub\ (expressed in cents of Euros, \euro), are reported in Table~\ref{table:results-random} in Appendix \ref{App:B}, together with the complete results for all methods.}
Regarding \matheur\ \rev{and SA}, \%gapLB is computed similarly to \milpmodel\ using the objective value of the best solution found during \rev{each algorithm's} execution, $z_{\textrm{ILS}}$\rev{\ (resp.\ $z_{\textrm{SA}}$)}, with respect to the \milpmodel\ lower bound $z_{\textrm{lb}}$ as $100\times ($$z_{\textrm{ILS}}$ - $z_{\textrm{lb}})/z_{\textrm{ILS}}$ (\rev{resp.\ using $z_{\textrm{SA}}$}). Finally, \inctime\ indicates the time required to reach the best solution found during the execution. \rev{For all three methods}, computation times are reported in seconds. 
\rev{A direct comparison of incumbent solutions against \milpmodel, expressed as \%gapUB for \matheur\ and SA, is provided in Table~\ref{table:results-random} in Appendix \ref{App:B}, together with \avglb\ and \avgub. This gap, computed as $100\times(z_{\textrm{ILS}} - z_{\textrm{ub}})/z_{\textrm{ILS}}$ (resp.\ using $z_{\textrm{SA}}$), is restricted to instances for which a feasible solution of MILP exists, and is omitted here for compactness.}

Table~\ref{table:results-random-summary} presents the results produced by the single execution of \milpmodel\ 
and the average of \rev{ten} executions of \matheur\ \rev{and SA} with different random seeds for the  variable consumption instances in the corresponding row. 
In this experiment, feasible solutions are not always available for \milpmodel, and a deeper inspection of its execution shows that for many instances it is not able to solve the initial relaxation of LP within the time limit (i.e., 1200 s). For this reason, the average \%gapLB is computed over a subset of 
instances within each group, and the number of instances considered is indicated in parentheses. Both \matheur\ and \rev{SA} show a more consistent behavior than \milpmodel, finding a feasible solution  every time the MILP does and identifying \rev{11.4} additional feasible solutions on average (in bold 
in column \#feas). \rev{In terms of computing time, \matheur\ is on average approximately 59\% faster than \milpmodel, and SA approximately 78\% faster; however, this shorter running time for SA comes at  the cost of solution quality: on average, \matheur\ improves on \milpmodel's incumbent by $2.06\%$ 
(\%gapUB $=-2.06\%$), while SA's incumbent solutions remain on average $10.23\%$ worse than \milpmodel's (\%gapUB $=+10.23\%$), as reported in Table~\ref{table:results-random} in Appendix \ref{App:B}.}
%
\begin{table}[ht]
\centering
\caption{Performance comparison of \milpmodel, \matheur, and \rev{SA} on instances with variable consumption. 
\comptime\ and \inctime\ are expressed in seconds. \matheur\ and \rev{SA} values \nfeas, \nopt, \%gapLB, 
and \inctime\ refer to the average results over \rev{10} executions. Numbers between parentheses for \%gapLB 
indicate the number of instances on which each value is computed in each execution. 
\rev{Full results, including \avglb, \avgub, and \%gapUB, are reported in 
Table~\ref{table:results-random} in Appendix \ref{App:B}.}}
\resizebox{\textwidth}{!}{%
\begin{tabular}{llc rrrr rrrr rrrr}
\toprule
\multicolumn{3}{c}{Instances} & \multicolumn{4}{c}{\milpmodel} & \multicolumn{4}{c}{\matheur} & \multicolumn{4}{c}{\rev{SA}} \\
\cmidrule(lr){1-3} \cmidrule(lr){4-7} \cmidrule(lr){8-11} \cmidrule(lr){12-15}
$|T|$ & $(|J|$, $|I|)$ & \ninst & \nfeas & \nopt & \%gapLB & \comptime & \nfeas & \nopt & \%gapLB & \inctime & \rev{\nfeas} & \rev{\nopt} & \rev{\%gapLB} & \rev{\inctime} \\
\cmidrule(lr){1-3} \cmidrule(lr){4-7} \cmidrule(lr){8-11} \cmidrule(lr){12-15}
\multirow{3}*{48} 
& ([5, 20], [3, 7]) & 9 & \textbf{7} & 4 & (7) 1.91 & 720.7 & \rev{\textbf{7.0}} & \rev{4.0} & \rev{(7) 1.88} & \rev{199.9} & \rev{\textbf{7.0}} & \rev{3.0} & \rev{(7) 11.98} & \rev{8.3} \\
& ([30, 90], [10, 20])  & 9 & \textbf{9} & 1 & (9) 1.06 & 1079.6 & \rev{\textbf{9.0}} & \rev{0.0} & \rev{(9) 1.24} & \rev{445.7} & \rev{\textbf{9.0}} & \rev{0.0} & \rev{(9) 14.01} & \rev{170.2} \\
& ([100, 200], [25, 35])  & 9 & \textbf{9} & 0 & (9) 0.26 & 1200.0 & \rev{\textbf{9.0}} & \rev{0.0} & \rev{(9) 0.68} & \rev{522.9} & \rev{\textbf{9.0}} & \rev{0.0} & \rev{(9) 11.91} & \rev{408.8} \\
\cmidrule(lr){2-2} \cmidrule(lr){3-3} \cmidrule(lr){4-7} \cmidrule(lr){8-11} \cmidrule(lr){12-15}
&  & 27 & \textbf{25} & 5 & (25) 1.01 & 1000.1 &
\rev{\textbf{25.0}} & \rev{4.0} & \rev{(25) 1.22} & \rev{404.7} & \rev{\textbf{25.0}} & \rev{3.0} & \rev{(25) 12.69} & \rev{210.8} \\
\cmidrule(lr){1-3} \cmidrule(lr){4-7} \cmidrule(lr){8-11} \cmidrule(lr){12-15}
\multirow{3}*{72}
& ([5, 20], [3, 7]) & 9 & \textbf{8} & 3 & (8) 2.55 & 802.1 & \rev{\textbf{8.0}} & \rev{3.0} & \rev{(8) 2.56} & \rev{215.3} & \rev{\textbf{8.0}} & \rev{1.5} & \rev{(8) 9.41} & \rev{44.8} \\
& ([30, 90], [10, 20]) & 9 & \textbf{9} & 0 & (9) 1.17 & 1200.0 & \rev{\textbf{9.0}} & \rev{0.0} & \rev{(9) 1.48} & \rev{490.4} & \rev{\textbf{9.0}} & \rev{0.0} & \rev{(9) 11.44} & \rev{124.9} \\
& ([100, 200], [25, 35]) & 9 & 7 & 0 & (6) 3.58 & 1200.0 & 
\rev{\textbf{9.0}} & \rev{0.0} & \rev{(6) 1.01} & \rev{570.8} & \rev{\textbf{9.0}} & \rev{0.0} & \rev{(6) 16.36} & \rev{431.9} \\
\cmidrule(lr){2-3} \cmidrule(lr){4-7} \cmidrule(lr){8-11} \cmidrule(lr){12-15}
&  & 27 & 24 & 3 & (23) 2.28 & 1067.4 & 
\rev{\textbf{26.0}} & \rev{3.0} & \rev{(23) 1.73} & \rev{433.6} & \rev{\textbf{26.0}} & \rev{1.5} & \rev{(23) 12.02} & \rev{206.5} \\
\cmidrule(lr){1-3} \cmidrule(lr){4-7} \cmidrule(lr){8-11} \cmidrule(lr){12-15}
\multirow{3}*{96} 
& ([5, 20], [3, 7]) & 9 & 6 & 3 & (6) 2.38 & 806.1 & 
\rev{\textbf{6.2}} & \rev{3.0} & \rev{(6) 2.59} & \rev{194.7} & \rev{\textbf{6.2}} & \rev{1.1} & \rev{(6) 9.50} & \rev{60.1} \\
& ([30, 90], [10, 20]) & 9 & 8 & 0 & (8) 1.54 & 1200.0 & 
\rev{\textbf{9.0}} & \rev{0.0} & \rev{(8) 1.75} & \rev{485.3} & \rev{\textbf{9.0}} & \rev{0.0} & \rev{(8) 13.80} & \rev{180.6} \\
& ([100, 200], [25, 35]) & 9 & 6 & 0 & (3) 7.57 & 1200.0 & 
\rev{\textbf{9.0}} & \rev{0.0} & \rev{(3) 1.56} & \rev{564.5} & \rev{\textbf{9.0}} & \rev{0.0} & \rev{(3) 20.34} & \rev{475.0} \\
\cmidrule(lr){2-3} \cmidrule(lr){4-7} \cmidrule(lr){8-11} \cmidrule(lr){12-15}
&  & 27 & 20 & 3 & (17) 2.90 &  1068.7 & 
\rev{\textbf{24.2}} & \rev{3.0} & \rev{(17) 2.01} & \rev{440.3} & \rev{\textbf{24.2}} & \rev{1.1} & \rev{(17) 13.43} & \rev{259.2} \\
\cmidrule(lr){1-3} \cmidrule(lr){4-7} \cmidrule(lr){8-11} \cmidrule(lr){12-15}
\multirow{3}*{120} 
& ([5, 20], [3, 7]) & 9 & \textbf{5} & 3 & (5) 1.71 & 848.6 & 
\rev{\textbf{5.0}} & \rev{3.0} & \rev{(5) 2.38} & \rev{246.9} & \rev{\textbf{5.0}} & \rev{2.0} & \rev{(5) 8.39} & \rev{120.0} \\
& ([30, 90], [10, 20]) & 9 & 5 & 0 & (5) 1.37 & 1200.0 & 
\rev{\textbf{6.7}} & \rev{0.0} & \rev{(5) 1.99} & \rev{516.5} & \rev{\textbf{6.7}} & \rev{0.0} & \rev{(5) 21.67} & \rev{226.8} \\
& ([100, 200], [25, 35]) & 9 & 3 & 0 & (1) 0.81 & 1200.0 & 
\rev{\textbf{6.6}} & \rev{0.0} & \rev{(1) 1.02} & \rev{573.9} & \rev{\textbf{6.5}} & \rev{0.0} & \rev{(1) 13.39} & \rev{479.4} \\
\cmidrule(lr){2-3} \cmidrule(lr){4-7} \cmidrule(lr){8-11} \cmidrule(lr){12-15}
&  & 27 & 13 & 3 & (11) 1.47 & 1082.9 & 
\rev{\textbf{18.2}} & \rev{3.0} & \rev{(11) 2.08} & \rev{462.9} & \rev{\textbf{18.2}} & \rev{2.0} & \rev{(11) 14.88} & \rev{287.7} \\
\cmidrule(lr){1-3} \cmidrule(lr){4-7} \cmidrule(lr){8-11} \cmidrule(lr){12-15}
Tot & & 108 & 82 & 14 & (76) 1.88 & 1054.8 & 
\rev{\textbf{93.4}} & \rev{13.0} & \rev{(76) 1.68} & \rev{433.3} & \rev{\textbf{93.4}} & \rev{7.6} & \rev{(76) 12.97} & \rev{237.1} \\
\bottomrule
\end{tabular}%
}
\label{table:results-random-summary}
\end{table}

From a broader perspective, these experiments provide some interesting insights regarding the \problemname. The structure of the consumption function has a direct impact on the performance of the methods. Tackling small-sized instances through \milpmodel\ combined with \texttt{Gurobi} appears to be a good alternative, producing good quality solutions in reasonable computation times, regardless of the consumption function assumed for the jobs (see Appendix \ref{App:A} for analogous results on instances with fixed consumption functions). However, when considering larger instances and a variable consumption, its performance degrades significantly. The results achieved by \matheur\ are rather consistent, finding near-optimal solutions and outperforming the results obtained by \milpmodel\ for larger instances (\%gapUB). In addition, these solutions show small percentage deviations from the best known \milpmodel\ lower bounds (\%gapLB), indicating high quality. 
 \rev{A similar comparison against SA further highlights the 
effectiveness of \matheur: despite SA's shorter running time, its solutions deviate on average $12.97\%$ from the \milpmodel\ lower bound, against $1.68\%$ for \matheur, and it only rarely improves on \milpmodel's incumbent solutions (\%gapUB $=+10.23\%$ on average, against $-2.06\%$ for \matheur). 
\matheur\ also reaches a proven optimal solution more often than SA (13.0 versus 7.6 instances on average, out of 108 instances).}

%
\rev{A statistical analysis of the difference between \matheur\ and SA on the small and medium instances is reported in Appendix \ref{App C} for compactness, and confirms that \matheur\ achieves significantly better average objective function values than SA. Overall, across our set of instances, \matheur\  outperforms both competing approaches: it is faster and finds more feasible solutions than \milpmodel, and it achieves better average objective function values than both \milpmodel\ and SA.
}

Figure \ref{fig:gaps} shows more in detail the variability of the \matheur\ gaps across the \rev{ten} executions averaged in Table \ref{table:results-random-summary} and Table \ref{table:results-random} of Appendix \ref{App:B}, specifically \%gapLB on the left (Figure \ref{fig:gap_lb}) and \%gapUB on the right (Figure \ref{fig:gap_ub}). Each group corresponds to a given 
time horizon and instance class (see size details in Table \ref{table:results-random-summary}), and shows the minimum, average, and maximum gap values observed. For each group, these statistics are computed on the subset of instances for which gaps exist in all \rev{ten} runs. Overall, \%gapLB values are variable but 
in a limited range between \rev{0.6\% and 3.0\%}. On the other hand, \%gapUB exhibits values close to 0 for most groups, at the exception of \rev{groups 72-3, 96-3, and 120-3}, showing \rev{negative} values, meaning that the heuristic is able to \rev{improve} the MILP upper bound provided within the time limit 
\rev{(most notably for 96-3 and 120-3, where the improvement exceeds 16\% on average)}.
\begin{figure}[ht]
\centering
\begin{subfigure}{0.49\textwidth}
\centering
\begin{tikzpicture}
\begin{axis}[
    width=\linewidth,
    height=6cm,
    ymin=0, ymax=3.5,
    xmin=0.5, xmax=12.5,
    ytick={0,0.5,1,1.5,2,2.5,3,3.5},
    yticklabel style={font=\scriptsize},
    xtick={1,2,3,4,5,6,7,8,9,10,11,12},
    xticklabels={
        48-1,48-2,48-3,
        72-1,72-2,72-3,
        96-1,96-2,96-3,
        120-1,120-2,120-3
    },
    x tick label style={
        rotate=90,
        anchor=east,
        font=\scriptsize
    },
    ymajorgrids=true,
    grid style={gray!30},
    axis lines=box
]
\addplot[black,no markers] coordinates {(1,1.75) (1,2.03)};
\addplot[black,no markers] coordinates {(2,1.18) (2,1.39)};
\addplot[black,no markers] coordinates {(3,0.63) (3,0.72)};
\addplot[black,no markers] coordinates {(4,2.32) (4,2.86)};
\addplot[black,no markers] coordinates {(5,1.36) (5,1.63)};
\addplot[black,no markers] coordinates {(6,0.97) (6,1.05)};
\addplot[black,no markers] coordinates {(7,2.20) (7,2.86)};
\addplot[black,no markers] coordinates {(8,1.58) (8,1.85)};
\addplot[black,no markers] coordinates {(9,1.44) (9,1.74)};
\addplot[black,no markers] coordinates {(10,2.09) (10,2.97)};
\addplot[black,no markers] coordinates {(11,1.81) (11,2.13)};
\addplot[black,no markers] coordinates {(12,0.89) (12,1.15)};
\addplot[
    black,
    only marks,
    mark=-,
    mark size=1.8pt,
    line width=0.6pt
] coordinates {
    (1,1.75) (2,1.18) (3,0.63)
    (4,2.32) (5,1.36) (6,0.97)
    (7,2.20) (8,1.58) (9,1.44)
    (10,2.09) (11,1.81) (12,0.89)
};
\addplot[
    black,
    only marks,
    mark=x,
    mark size=1.8pt,
    line width=0.6pt
] coordinates {
    (1,1.88) (2,1.24) (3,0.68)
    (4,2.56) (5,1.48) (6,1.01)
    (7,2.59) (8,1.75) (9,1.56)
    (10,2.38) (11,1.99) (12,1.02)
};
\addplot[
    black,
    only marks,
    mark=-,
    mark size=1.8pt,
    line width=0.6pt
] coordinates {
    (1,2.03) (2,1.39) (3,0.72)
    (4,2.86) (5,1.63) (6,1.05)
    (7,2.86) (8,1.85) (9,1.74)
    (10,2.97) (11,2.13) (12,1.15)
};
\end{axis}
\end{tikzpicture}
\caption{\%gapLB}
\label{fig:gap_lb}
\end{subfigure}
\hfill
\begin{subfigure}{0.49\textwidth}
\centering
\begin{tikzpicture}
\begin{axis}[
    width=\linewidth,
    height=6cm,
    ymin=-18, ymax=2,
    xmin=0.5, xmax=12.5,
    ytick={-18,-16,-14,-12,-10,-8,-6,-4,-2,0,2},
    yticklabel style={font=\scriptsize},
    xtick={1,2,3,4,5,6,7,8,9,10,11,12},
    xticklabels={
        48-1,48-2,48-3,
        72-1,72-2,72-3,
        96-1,96-2,96-3,
        120-1,120-2,120-3
    },
    x tick label style={
        rotate=90,
        anchor=east,
        font=\scriptsize
    },
    ymajorgrids=true,
    grid style={gray!30},
    axis lines=box
]
\addplot[black,no markers] coordinates {(1,-0.18) (1,0.12)};
\addplot[black,no markers] coordinates {(2,0.12) (2,0.34)};
\addplot[black,no markers] coordinates {(3,0.37) (3,0.46)};
\addplot[black,no markers] coordinates {(4,-0.25) (4,0.31)};
\addplot[black,no markers] coordinates {(5,0.19) (5,0.47)};
\addplot[black,no markers] coordinates {(6,-5.18) (6,-5.12)};
\addplot[black,no markers] coordinates {(7,-0.20) (7,0.50)};
\addplot[black,no markers] coordinates {(8,0.04) (8,0.32)};
\addplot[black,no markers] coordinates {(9,-16.99) (9,-16.81)};
\addplot[black,no markers] coordinates {(10,0.40) (10,1.33)};
\addplot[black,no markers] coordinates {(11,0.45) (11,0.79)};
\addplot[black,no markers] coordinates {(12,-16.42) (12,-16.26)};
\addplot[
    black,
    only marks,
    mark=-,
    mark size=1.8pt,
    line width=0.6pt
] coordinates {
    (1,-0.18) (2,0.12) (3,0.37)
    (4,-0.25) (5,0.19) (6,-5.18)
    (7,-0.20) (8,0.04) (9,-16.99)
    (10,0.40) (11,0.45) (12,-16.42)
};
\addplot[
    black,
    only marks,
    mark=x,
    mark size=1.8pt,
    line width=0.6pt
] coordinates {
    (1,-0.03) (2,0.19) (3,0.42)
    (4,0.01) (5,0.31) (6,-5.15)
    (7,0.22) (8,0.21) (9,-16.90)
    (10,0.71) (11,0.63) (12,-16.33)
};
\addplot[
    black,
    only marks,
    mark=-,
    mark size=1.8pt,
    line width=0.6pt
] coordinates {
    (1,0.12) (2,0.34) (3,0.46)
    (4,0.31) (5,0.47) (6,-5.12)
    (7,0.50) (8,0.32) (9,-16.81)
    (10,1.33) (11,0.79) (12,-16.26)
};
\end{axis}
\end{tikzpicture}
\caption{\%gapUB}
\label{fig:gap_ub}
\end{subfigure}
\caption{Variation ranges of \matheur\ gaps over \rev{ten} executions. Minimum and maximum values for each group are indicated by bars, while crosses indicate means.}
\label{fig:gaps}
\end{figure}

\subsubsection{Evaluation of \matheur\ on large instances}
\label{subsec:res-larger}
To evaluate the performance of the proposed \matheur\ on large-sized instances, we consider an adapted version of the algorithm described in Algorithm~\ref{alg:ils}, in which the application of the \milpmodel\ within the constructive phase (Steps~2 and~3 in Algorithm~\ref{alg:insertion}, Section~\ref{sec:heur:constructive}) is omitted. The experiments are conducted over the 45 large instances with variable consumption, defined by $(|J|,|I|,|T|) = (300,45,144)$. The test set is partitioned into three groups of 15 instances each, corresponding to saturation levels $\nu \in {0.7,0.8,0.9}$, and again for each instance we impose a time limit of 600 seconds. As discussed in the previous section, for large instances the \milpmodel\ is not effective due to the size of the resulting formulation.

As a benchmark, we compare our approach with \rev{the same SA reported in Section \ref{sec:overall-comparison}, but which
}
%
%
relies on the same constructive scheme as the simplified \matheur\ considered in this experiment. 

%
%
%
\begin{table}[ht]
\centering
\caption{Performance comparison between SA and \matheur\ on large instances with variable energy consumption, characterized by $(|J|,|I|,|T|) = (300,45,144)$ and saturation levels $\nu \in \{0.7,0.8,0.9\}$. Values of Avg. Best are expressed in cents of \euro. Values of \inctime\ are expressed in seconds. Numbers between parentheses for \%gapBest indicate the number of instances on which each value is computed in each of the \rev{10} runs of execution.}
\setlength{\tabcolsep}{6pt}
\renewcommand{\arraystretch}{0.6}
\resizebox{\textwidth}{!}{%

\begin{tabular}{clc rrrr rrrr}
\toprule
\multicolumn{3}{c}{\small Instances} & \multicolumn{4}{c}{\small SA} & \multicolumn{4}{c}{\small ILS} \\
\cmidrule(lr){1-3} \cmidrule(lr){4-7} \cmidrule(lr){8-11}
{\small Run} & {\small $\nu$} & {\small \#inst} & {\small \#feas} & {\small Avg. Best }& {\small inc.time} & {\small \rev{RPD (\%)}} & {\small \#feas}  & {\small inc.time} & {\small \rev{RPD (\%)}} & {\small \rev{\%gapSA}} \\
\cmidrule(lr){1-3} \cmidrule(lr){4-7} \cmidrule(lr){8-11}
\multirow{3}*{{\small Avg}}
& {\small 0.7} & {\small 15} & {\small \rev{15}} & {\small \rev{87672.72}} & {\small \rev{597.9}} & {\small \rev{23.3}}  & {\small \rev{15}} & {\small \rev{401.1}} & {\small \rev{0.05}} & {\small \rev{(15) -23.2$^{***}$}} \\
& {\small 0.8} & {\small 15} & {\small \rev{15}} & {\small \rev{100014.49}} & {\small \rev{597.4}} & {\small \rev{16.4}}  & {\small \rev{15}} & {\small \rev{370.6}} & {\small \rev{0.08}} & {\small \rev{(15) -16.3$^{***}$}} \\
& {\small 0.9} & {\small 15} & {\small \rev{15}} & {\small \rev{114951.46}} & {\small \rev{594.5}} & {\small \rev{4.9}}  & {\small \rev{15}} & {\small \rev{227.8}} & {\small \rev{0.20}} & {\small \rev{(15) \phantom{0}-4.7$^{***}$}}  \\
\cmidrule(lr){1-3} \cmidrule(lr){4-7} \cmidrule(lr){8-11}
{\small Tot} & & {\small 45} & {\small \rev{45}} & {\small \rev{100879.56}} & {\small \rev{596.6}} & {\small \rev{14.9}} & {\small \rev{45}} & {\small \rev{333.2}} & {\small \rev{0.11}} & {\small \rev{(45) -14.7$^{***}$}}  \\
\bottomrule
\end{tabular}%
}
\label{tab:large_instances_runs}
\end{table}
Table~\ref{tab:large_instances_runs} reports the average results of \rev{ten} independent runs of both the simplified \matheur\ and the SA algorithm, grouped by the saturation level $\nu$. 
The main message is that \matheur\ consistently produces higher-quality solutions, with average gaps from the average best solutions found by SA \rev{(column \%gapSA)} of up to \rev{$-23.2\%$, $-16.3\%$, and $-4.7\%$} as the saturation level increases, respectively. 
\rev{To assess the statistical significance of this difference, we compare, for each of the 45 instances, the average objective value obtained by \matheur\ and SA over the ten runs. A Shapiro-Wilk test rejects normality of the paired per-instance differences ($p = 0.04$), so we apply the non-parametric Wilcoxon signed-rank test, following \citet{Osaba2021}. The test confirms that \matheur\ significantly outperforms SA ($p < 0.001$, marked by $^{***}$ in column \%gapSA of Table~\ref{tab:large_instances_runs}); in fact, \matheur\ achieves a better average objective value than SA on all 45 instances.}
\rev{We also report, for each method, the relative percentage difference (RPD) between its average objective value over the ten runs ($\bar{z}_{\text{method}} $) and the best known value for the instance across both methods ($z_{\text{BKS}}$), i.e., $100\times(\bar{z}_{\text{method}} - z_{\text{BKS}})/z_{\text{BKS}}$. The RPD of \matheur\ remains below $0.2\%$ in every group (overall $0.11\%$), confirming that its average performance is essentially indistinguishable from its own best solution, while the RPD of SA ranges between $4.9\%$ and $23.3\%$ (overall $14.9\%$), reflecting a substantial and consistent gap from the best known solutions, which are found by \matheur\ on all 45 instances.}

\subsection{Impact of variable consumption }
\label{subsec:qualitative-results}
We now focus on the specific impact on feasibility and budget constraints \eqref{eq6_s4} when incorporating variable consumptions on our generated instances.
Given a base configuration as described in Section~\ref{subsec:instances}, consider the best known solution $s$ for the fixed consumption. Solution $s$ specifies the machine assignments and the starting times for each job.
Solution $s$ is also feasible in terms of pure scheduling (as the processing times remain the same) and fits within the planning horizon when transferred to the variable consumption setting using the same base configuration. However, it may be infeasible with respect to the energy budget due to the variable energy profiles of jobs.
Therefore, for every tested base configuration, we evaluate different feasibility metrics when plugging a solution obtained for a fixed consumption pattern into the corresponding variable-consumption instance. In other words, we fix the starting times and machine assignments from the fixed-consumption solution and apply them directly to each related variable-consumption instance.
This allows us to measure the impact of introducing variable job consumption profiles on the feasibility of solutions originally designed under fixed-consumption assumptions, and evaluate how realistic variations in energy usage affect the feasibility and energy efficiency of scheduling solutions.

\newcommand{\ninf}{\#inf}
\newcommand{\mininf}{min inf}
\newcommand{\avginf}{avg inf}
\newcommand{\maxinf}{max inf}
\newcommand{\totalinf}{total inf}
\newcommand{\percinc}{\%infE} 
\newcommand{\percinf}{\%infT} 

Table~\ref{tab:quality_comparison} and Figure \ref{fig:quality_comparison} present the aggregated results of one execution of this analysis on small-to-medium size instances. In Table~\ref{tab:quality_comparison}, \mininf, \avginf, and \maxinf\ are the average per group of minimum, average, and maximum amount of energy exceeding the budget; \totalinf\ stands for the average per group of the total sum of energy above this budget. In this experiment, groups are identified by the time horizon $|T|$.
For a given instance, let $\bar{E}_t$ be the amount of energy scheduled in time slot $t$. We set \percinc\ $= \sum\limits_{t \in T}\max(0, \bar{E}_t-E)/\sum\limits_{t \in T}\bar{E}_t$ as the percentage by which the energy budget $E$ needs to be increased for the instance to become feasible, while \percinf\ $= |\{t \in T: \bar{E}_t > E\}|/|T|$ is the percentage of infeasible time slots. 
In Figure~\ref{fig:quality_comparison}, we show the average \percinc\ and \percinf\ over all groups (with respect to $|T|$) and overall (Tot).

None of the solutions obtained for the instances with fixed consumption is feasible when evaluated in the corresponding variable consumption instance, showing the challenges of our problem with respect to the state-of-the-art parallel machine scheduling problems with constant resource consumption functions. 
The number of infeasible time slots, on average, represents less than 20\% of the planning horizon. However, the increase in energy budget required to achieve feasibility is significant: approximately 35\% for variable instances. 
This highlights the importance of considering variable consumption functions in energy-aware scheduling studies. 
In real-world applications, ignoring such variations would require significantly increasing the energy budget to regain feasibility or, alternatively, removing some jobs from the schedule to avoid exceeding the maximum budget. In general, neglecting variable consumption at the planning stage has a direct impact on the quality of the schedule and its potential practical use.

\begin{minipage}{0.5\textwidth}
    \centering
    \scriptsize
    \centering
    \captionof{table}{Energy exceeding the budget expressed in kWh.}
    \resizebox{0.8\textwidth}{!}{
    \begin{tabular}{cccccc}
    \toprule
    \# & \(|T|\) & min inf & avg inf & max inf & total inf \\
    \cmidrule(lr){1-2} \cmidrule(lr){3-6}
    27 & 48  & 15   & 104   & 296  & 946  \\
    27 & 72  & 10   & 99   & 263   & 1190  \\
    27 & 96  & 9   & 93   & 274   & 1561  \\
    27 & 120 & 4   & 78   & 243   & 2015  \\
    \cmidrule(lr){1-2} \cmidrule(lr){3-6}
    Tot &  & 10   & 93   & 269   & 1428  \\
    \bottomrule
    \end{tabular}
    }
    \label{tab:quality_comparison}
\end{minipage}
\hfill
\begin{minipage}{0.5\textwidth}
    \centering
    \resizebox{\textwidth}{!}{ 
        \begin{tikzpicture}
            \draw[-] (0,0) -- (10,0) node[anchor=west] {}; 
            \draw[-] (0,0) -- (0,5.3) node[anchor=south] {};       

            \node[anchor=east] at (-0.2,0.5) {48};
            \node[anchor=east] at (-0.2,1.5) {72};
            \node[anchor=east] at (-0.2,2.5) {96};
            \node[anchor=east] at (-0.2,3.5) {120};
            \node[anchor=east] at (-0.2,4.5) {Tot};

            \fill[gray!20] (0,0.4) rectangle (7.8,0.6) node[midway,below,black] {\footnotesize 39\%};
            \fill[gray!20] (0,1.4) rectangle (7,1.6) node[midway,below,black] {\footnotesize 35\%};
            \fill[gray!20] (0,2.4) rectangle (7.4,2.6) node[midway,below,black] {\footnotesize 37\%};
            \fill[gray!20] (0,3.4) rectangle (6.4,3.6) node[midway,below,black] {\footnotesize 32\%};
            \fill[gray!20] (0,4.4) rectangle (7.2,4.6) node[midway,below,black] {\footnotesize 36\%};

            \fill[gray!40] (0,0.6) rectangle (3.8,0.8) node[midway,above,black] {\footnotesize 19\%};
            \fill[gray!40] (0,1.6) rectangle (3.2,1.8) node[midway,above,black] {\footnotesize 16\%};
            \fill[gray!40] (0,2.6) rectangle (3.4,2.8) node[midway,above,black] {\footnotesize 17\%};
            \fill[gray!40] (0,3.6) rectangle (4,3.8) node[midway,above,black] {\footnotesize 20\%};
            \fill[gray!40] (0,4.6) rectangle (3.6,4.8) node[midway,above,black] {\footnotesize 18\%};

            \foreach \x/\name in {0/0, 1/5, 2/10, 3/15, 4/20, 5/25, 6/30, 7/35, 8/40, 9/45, 10/50} {
                \draw[shift={(\x,0)}] (0,0.0) -- (0,-0.2) node[anchor=north] {\footnotesize \name\%};
            }

            \draw[fill=gray!40] (9.0,4.7) rectangle (9.5,4.9);
            \node[anchor=west] at (9.55,4.8) {\%infT};

            \draw[fill=gray!20] (9.0,4.5) rectangle (9.5,4.3);
            \node[anchor=west] at (9.55,4.4) {\%infE};
        \end{tikzpicture}
    }
    
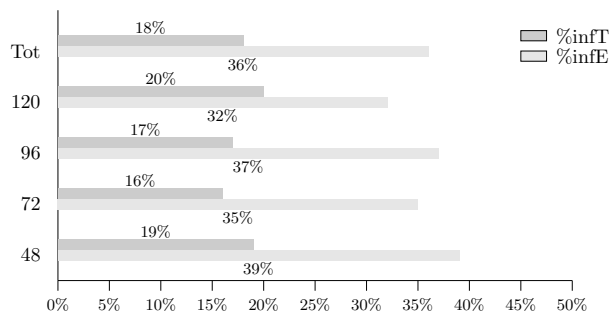
\captionof{figure}{Percentage of infeasibility.}
    \label{fig:quality_comparison}
\end{minipage}%

\begin{figure}
    \centering
    \begin{subfigure}[t]{0.45\textwidth}
        \centering
        \resizebox{\textwidth}{!}{ 
        \begin{tikzpicture}[]    
        \draw[->, line width=10pt] (0,0) -- (95,0) node[right] {\scalebox{7}{\textbf{T}}}; 
        \draw[->, line width=10pt] (0,0) -- (0,37.5) node[left, xshift=-30pt] {\scalebox{7}{\textbf{kWh}}}; 

        \foreach \y/\name in {7.5/120, 15/240, 22.5/360, 30/480} 
            \draw (-2,\y) node[left] {\scalebox{7}{\name}}; 
        
        \foreach \x/\name in {0/0:00, 28/7:00, 72/18:00,  95/24:00} 
            \draw (\x,-2) node[left] {\scalebox{7}{\name}}; 

        \draw[thick, gray!50, dash pattern=on 70pt off 25pt, line width=5] (0,23.125) -- (95,23.125) node[right] {};
        
        \draw[thick, black]
        (0,23.125) -- (1,23.125)--(1,23.125) -- (2,23.125)--(2,23.125) -- (3,23.125)--(3,23.125) -- (4,23.125)--(4,23.125) -- (5,23.125)--(5,23.125) -- (6,23.125)--(6,23.125) -- (7,23.125)--(7,23.125) -- (8,23.125)--(8,23.125) -- (9,23.125)--(9,23.125) -- (10,23.125)--(10,23.125) -- (11,23.125)--(11,23.125) -- (12,23.125)--(12,23.125) -- (13,23.125)--(13,23.125) -- (14,23.125)--(14,23.125) -- (15,23.125)--(15,23.125) -- (16,23.125)--(16,23.125) -- (17,23.125)--(17,23.125) -- (18,23.125)--(18,23.125) -- (19,23.125)--(19,23.125) -- (20,23.125)--(20,23.125) -- (21,23.125)--(21,23.125) -- (22,23.125)--(22,23.125) -- (23,23.125)--(23,23.125) -- (24,23.125)--(24,23.125) -- (25,23.125)--(25,23.125) -- (26,23.125)--(26,23.125) -- (27,23.125)--(27,23.125) -- (28,23.125)--(28,23.125) -- (29,23.125)--(29,23.125) -- (30,23.125)--(30,23.125) -- (31,23.125)--(31,23.125) -- (32,23.125)--(32,18.75) -- (33,18.75)--(33,14.375) -- (34,14.375)--(34,14.375) -- (35,14.375)--(35,14.375) -- (36,14.375)--(36,14.375) -- (37,14.375)--(37,14.375) -- (38,14.375)--(38,11.25) -- (39,11.25)--(39,8.125) -- (40,8.125)--(40,8.125) -- (41,8.125)--(41,5.0) -- (42,5.0)--(42,5.0) -- (43,5.0)--(43,8.125) -- (44,8.125)--(44,8.125) -- (45,8.125)--(45,8.125) -- (46,8.125)--(46,11.25) -- (47,11.25)--(47,11.25) -- (48,11.25)--(48,11.25) -- (49,11.25)--(49,11.25) -- (50,11.25)--(50,14.375) -- (51,14.375)--(51,14.375) -- (52,14.375)--(52,14.375) -- (53,14.375)--(53,14.375) -- (54,14.375)--(54,14.375) -- (55,14.375)--(55,14.375) -- (56,14.375)--(56,14.375) -- (57,14.375)--(57,14.375) -- (58,14.375)--(58,14.375) -- (59,14.375)--(59,14.375) -- (60,14.375)--(60,14.375) -- (61,14.375)--(61,14.375) -- (62,14.375)--(62,14.375) -- (63,14.375)--(63,14.375) -- (64,14.375)--(64,14.375) -- (65,14.375)--(65,11.25) -- (66,11.25)--(66,11.25) -- (67,11.25)--(67,11.25) -- (68,11.25)--(68,8.125) -- (69,8.125)--(69,8.125) -- (70,8.125)--(70,8.125) -- (71,8.125)--(71,8.125) -- (72,8.125)--(72,8.125) -- (73,8.125)--(73,5.0) -- (74,5.0)--(74,5.0) -- (75,5.0)--(75,5.0) -- (76,5.0)--(76,1.875) -- (77,1.875)--(77,1.875) -- (78,1.875)--(78,8.125) -- (79,8.125)--(79,8.125) -- (80,8.125)--(80,11.25) -- (81,11.25)--(81,14.375) -- (82,14.375)--(82,14.375) -- (83,14.375)--(83,14.375) -- (84,14.375)--(84,14.375) -- (85,14.375)--(85,14.375) -- (86,14.375)--(86,14.375) -- (87,14.375)--(87,14.375) -- (88,14.375)--(88,14.375) -- (89,14.375)--(89,14.375) -- (90,14.375)--(90,14.375) -- (91,14.375)--(91,14.375) -- (92,14.375)--(92,14.375) -- (93,14.375)--(93,14.375) -- (94,14.375)--(94,14.375) -- (95,14.375)
        node[right] {};
        
        \draw[thick, gray!80,  dash pattern=on 28pt off 10pt, line width = 5pt]
        (0,0) -- (28,0) -- (29,0) -- (29,0.8125) -- (30,0.8125) -- (30,1.625) -- (31,1.625) -- (31,2.4375) -- (32,2.4375) -- (32,3.25) -- (33,3.25) -- (33,4.0625) -- (34,4.0625) -- (34,4.875) -- (35,4.875) -- (35,5.6875) -- (36,5.6875) -- (36,6.5) -- (37,6.5) -- (37,7.3125) -- (38,7.3125) -- (38,8.125) -- (39,8.125) -- (39,8.9375) -- (40,8.9375) -- (40,9.75) -- (41,9.75) -- (41,10.5625) -- (42,10.5625) -- (42,11.375) -- (43,11.375) -- (43,12.1875) -- (44,12.1875) -- (44,13.0) -- (45,13.0) -- (45,13.8125) -- (46,13.8125) -- (46,14.625) -- (47,14.625) -- (47,15.4375) -- (48,15.4375) -- (48,16.25) -- (49,16.25) -- (49,17.0625) -- (50,17.0625) -- (50,16.25) -- (51,16.25) -- (51,15.4375) -- (52,15.4375) -- (52,14.625) -- (53,14.625) -- (53,13.8125) -- (54,13.8125) -- (54,13.0) -- (55,13.0) -- (55,12.1875) -- (56,12.1875) -- (56,11.375) -- (57,11.375) -- (57,10.5625) -- (58,10.5625) -- (58,9.75) -- (59,9.75) -- (59,8.9375) -- (60,8.9375) -- (60,8.125) -- (61,8.125) -- (61,7.3125) -- (62,7.3125) -- (62,6.5) -- (63,6.5) -- (63,5.6875) -- (64,5.6875) -- (64,4.875) -- (65,4.875) -- (65,4.0625) -- (66,4.0625) -- (66,3.25) -- (67,3.25) -- (67,2.5625) -- (68,2.5625) -- (68,1.625) -- (69,1.625) -- (69,0.8125) -- (70,0.8125) -- (70,0.0) -- (95,0.0)
        node[right] {};
        
        \draw[gray!50, dash pattern=on 70pt off 25pt, line width=5pt] (42.5,35) -- ++(5,0) node[right] {\scalebox{7}{\textbf{Energy budget}}};
        \draw[black] (42.5,40) -- ++(5,0) node[right] {\scalebox{7}{\textbf{Consumption}}};
        \draw[gray!80, dash pattern=on 28pt off 10pt, line width = 5pt] (42.5,30) -- ++(5,0) node[right] {\scalebox{7}{\textbf{Photovoltaic energy}}};
        
        \end{tikzpicture}
        }
        \caption{Optimal schedule for fixed consumptions.}
        \label{fig:6a}
    \end{subfigure}
    \vfill
    \vspace{.5cm}
    \begin{subfigure}[t]{0.45\textwidth}
        \centering
        \resizebox{\textwidth}{!}{ 
        \begin{tikzpicture}[]    
        \draw[->, line width=10pt] (0,0) -- (95,0) node[right] {\scalebox{7}{\textbf{T}}}; 
        \draw[->, line width=10pt] (0,0) -- (0,37.5) node[left, xshift=-30pt] {\scalebox{7}{\textbf{kWh}}}; 

        \foreach \y/\name in {7.5/120, 15/240, 22.5/360, 30/480} 
            \draw (-2,\y) node[left] {\scalebox{7}{\name}}; 
        
        \foreach \x/\name in {0/0:00, 28/7:00, 72/18:00,  95/24:00} 
            \draw (\x,-2) node[left] {\scalebox{7}{\name}}; 

        \draw[thick, gray!50, dash pattern=on 70pt off 25pt, line width=5] (0,23.125) -- (95,23.125) node[right] {};
        
        \draw[thick, black]
        (0,36.75) -- (1,36.75) -- (1,29.25) -- (2,29.25) -- (2,22.25) -- (3,22.25)
        -- (3,26.0625) -- (4,26.0625) -- (4,32.1875) -- (5,32.1875) -- (5,22.25)
        -- (6,22.25) -- (6,22.3125) -- (7,22.3125) -- (7,13.6875) -- (8,13.6875)
        -- (8,25.125) -- (9,25.125) -- (9,21.5) -- (10,21.5) -- (10,34.75)
        -- (11,34.75) -- (11,26.4375) -- (12,26.4375) -- (12,28.1875) -- (13,28.1875)
        -- (13,26.6875) -- (14,26.6875) -- (14,29.375) -- (15,29.375) -- (15,27.125)
        -- (16,27.125) -- (16,25.25) -- (17,25.25) -- (17,11.25) -- (18,11.25)
        -- (18,11.9375) -- (19,11.9375) -- (19,26.3125) -- (20,26.3125) -- (20,17.3125)
        -- (21,17.3125) -- (21,19.1875) -- (22,19.1875) -- (22,12.3125) -- (23,12.3125)
        -- (23,22.75) -- (24,22.75) -- (24,25.5) -- (25,25.5) -- (25,18.875)
        -- (26,18.875) -- (26,15.25) -- (27,15.25) -- (27,24.4375) -- (28,24.4375)
        -- (28,12.4375) -- (29,12.4375) -- (29,13.25) -- (30,13.25) -- (30,6.0625)
        -- (31,6.0625) -- (31,10.5) -- (32,10.5) -- (32,9.125) -- (33,9.125)
        -- (33,12.625) -- (34,12.625) -- (34,9.4375) -- (35,9.4375) -- (35,7.25)
        -- (36,7.25) -- (36,3.5625) -- (37,3.5625) -- (37,6.25) -- (38,6.25)
        -- (38,10.5) -- (39,10.5) -- (39,5.4375) -- (40,5.4375) -- (40,12.6875)
        -- (41,12.6875) -- (41,9.625) -- (42,9.625) -- (42,4.875) -- (43,4.875)
        -- (43,4.8125) -- (44,4.8125) -- (44,18.75) -- (45,18.75) -- (45,18.0625)
        -- (46,18.0625) -- (46,10.875) -- (47,10.875) -- (47,15.75) -- (48,15.75)
        -- (48,16.25) -- (49,16.25) -- (49,19.625) -- (50,19.625) -- (50,10.125)
        -- (51,10.125) -- (51,10.75) -- (52,10.75) -- (52,17.375) -- (53,17.375)
        -- (53,18.375) -- (54,18.375) -- (54,9.8125) -- (55,9.8125) -- (55,19.875)
        -- (56,19.875) -- (56,14.375) -- (57,14.375) -- (57,18.875) -- (58,18.875)
        -- (58,24.4375) -- (59,24.4375) -- (59,15.875) -- (60,15.875) -- (60,4.25)
        -- (61,4.25) -- (61,9.625) -- (62,9.625) -- (62,17.4375) -- (63,17.4375)
        -- (63,2.6875) -- (64,2.6875) -- (64,2.1875) -- (65,2.1875) -- (65,9.1875)
        -- (66,9.1875) -- (66,4.1875) -- (67,4.1875) -- (67,7.5625) -- (68,7.5625)
        -- (68,6.8125) -- (69,6.8125) -- (69,8.6875) -- (70,8.6875) -- (70,2.625)
        -- (71,2.625) -- (71,1.75) -- (72,1.75) -- (72,2.1875) -- (73,2.1875) -- (73,7.9375)
        -- (74,7.9375) -- (74,8.125) -- (75,8.125) -- (75,12.875) -- (76,12.875)
        -- (76,13.8125) -- (77,13.8125) -- (77,12.8125) -- (78,12.8125) -- (78,15.75)
        -- (79,15.75) -- (79,12.125) -- (80,12.125) -- (80,9.3125) -- (81,9.3125)
        -- (81,15.5) -- (82,15.5) -- (82,15.875) -- (83,15.875) -- (83,20.125)
        -- (84,20.125) -- (84,17.25) -- (85,17.25) -- (85,13.0) -- (86,13.0)
        -- (86,10.5) -- (87,10.5) -- (87,14.25) -- (88,14.25) -- (88,14.9375)
        -- (89,14.9375) -- (89,20.25) -- (90,20.25) -- (90,15.9375) -- (91,15.9375)
        -- (91,9.25) -- (92,9.25) -- (92,18.375) -- (93,18.375) -- (94,18.375)-- (94,7.375) -- (95,7.375) -- (95,7.375)
        node[right] {};
        
        \draw[thick, gray!80,  dash pattern=on 28pt off 10pt, line width = 5pt]
        (0,0) -- (28,0) -- (29,0) -- (29,0.8125) -- (30,0.8125) -- (30,1.625) -- (31,1.625) -- (31,2.4375) -- (32,2.4375) -- (32,3.25) -- (33,3.25) -- (33,4.0625) -- (34,4.0625) -- (34,4.875) -- (35,4.875) -- (35,5.6875) -- (36,5.6875) -- (36,6.5) -- (37,6.5) -- (37,7.3125) -- (38,7.3125) -- (38,8.125) -- (39,8.125) -- (39,8.9375) -- (40,8.9375) -- (40,9.75) -- (41,9.75) -- (41,10.5625) -- (42,10.5625) -- (42,11.375) -- (43,11.375) -- (43,12.1875) -- (44,12.1875) -- (44,13.0) -- (45,13.0) -- (45,13.8125) -- (46,13.8125) -- (46,14.625) -- (47,14.625) -- (47,15.4375) -- (48,15.4375) -- (48,16.25) -- (49,16.25) -- (49,17.0625) -- (50,17.0625) -- (50,16.25) -- (51,16.25) -- (51,15.4375) -- (52,15.4375) -- (52,14.625) -- (53,14.625) -- (53,13.8125) -- (54,13.8125) -- (54,13.0) -- (55,13.0) -- (55,12.1875) -- (56,12.1875) -- (56,11.375) -- (57,11.375) -- (57,10.5625) -- (58,10.5625) -- (58,9.75) -- (59,9.75) -- (59,8.9375) -- (60,8.9375) -- (60,8.125) -- (61,8.125) -- (61,7.3125) -- (62,7.3125) -- (62,6.5) -- (63,6.5) -- (63,5.6875) -- (64,5.6875) -- (64,4.875) -- (65,4.875) -- (65,4.0625) -- (66,4.0625) -- (66,3.25) -- (67,3.25) -- (67,2.5625) -- (68,2.5625) -- (68,1.625) -- (69,1.625) -- (69,0.8125) -- (70,0.8125) -- (70,0.0) -- (95,0.0)
        node[right] {};
        
        \draw[gray!50, dash pattern=on 70pt off 25pt, line width=5pt] (42.5,35) -- ++(5,0) node[right] {\scalebox{7}{\textbf{Energy budget}}};
        \draw[black] (42.5,40) -- ++(5,0) node[right] {\scalebox{7}{\textbf{Consumption}}};
        \draw[gray!80, dash pattern=on 28pt off 10pt, line width = 5pt] (42.5,30) -- ++(5,0) node[right] {\scalebox{7}{\textbf{Photovoltaic energy}}};
        
        \end{tikzpicture}
        }
        \caption{Optimal schedule for fixed consumptions re-evaluated on variable consumptions.}
        \label{fig:6b}
    \end{subfigure}
    \hfill
    \begin{subfigure}[t]{0.45\textwidth}
        \centering
        \resizebox{\textwidth}{!}{ 
        \begin{tikzpicture}[]

\draw[->, line width=10pt] (0,0) -- (95,0) node[right] {\scalebox{7}{\textbf{T}}}; 
\draw[->, line width=10pt] (0,0) -- (0,37.5) node[left, xshift=-30pt] {\scalebox{7}{\textbf{kWh}}}; 


\foreach \y/\name in {7.5/120, 15/240, 22.5/360, 30/480} 
    \draw (-2,\y) node[left] {\scalebox{7}{\name}}; 

\foreach \x/\name in {0/0:00, 28/7:00, 72/18:00,  95/24:00} 
    \draw (\x,-2) node[left] {\scalebox{7}{\name}}; 

\draw[thick, gray!50, dash pattern=on 70pt off 25pt, line width=5] (0,23.125) -- (95,23.125) node[right] {};

\draw[thick, black] 
(0, 21.875) -- (1, 21.875)--(1, 13.8125) -- (2, 13.8125)--(2, 11.8125) -- (3, 11.8125)--(3, 13.6875) -- (4, 13.6875)--(4, 13.3125) -- (5, 13.3125)--(5, 21.75) -- (6, 21.75)--(6, 20.125) -- (7, 20.125)--(7, 9.5) -- (8, 9.5)--(8, 14.0625) -- (9, 14.0625)--(9, 22.0) -- (10, 22.0)--(10, 22.3125) -- (11, 22.3125)--(11, 16.125) -- (12, 16.125)--(12, 18.1875) -- (13, 18.1875)--(13, 22.5625) -- (14, 22.5625)--(14, 22.0625) -- (15, 22.0625)--(15, 21.0) -- (16, 21.0)--(16, 19.75) -- (17, 19.75)--(17, 22.25) -- (18, 22.25)--(18, 13.9375) -- (19, 13.9375)--(19, 15.8125) -- (20, 15.8125)--(20, 14.9375) -- (21, 14.9375)--(21, 11.125) -- (22, 11.125)--(22, 20.3125) -- (23, 20.3125)--(23, 13.6875) -- (24, 13.6875)--(24, 18.5625) -- (25, 18.5625)--(25, 19.4375) -- (26, 19.4375)--(26, 14.875) -- (27, 14.875)--(27, 17.4375) -- (28, 17.4375)--(28, 11.9375) -- (29, 11.9375)--(29, 15.5625) -- (30, 15.5625)--(30, 13.8125) -- (31, 13.8125)--(31, 16.0625) -- (32, 16.0625)--(32, 6.0) -- (33, 6.0)--(33, 17.5625) -- (34, 17.5625)--(34, 6.625) -- (35, 6.625)--(35, 12.625) -- (36, 12.625)--(36, 9.125) -- (37, 9.125)--(37, 9.5) -- (38, 9.5)--(38, 9.375) -- (39, 9.375)--(39, 10.9375) -- (40, 10.9375)--(40, 11.0) -- (41, 11.0)--(41, 12.0) -- (42, 12.0)--(42, 10.1875) -- (43, 10.1875)--(43, 10.1875) -- (44, 10.1875)--(44, 12.0625) -- (45, 12.0625)--(45, 13.3125) -- (46, 13.3125)--(46, 17.9375) -- (47, 17.9375)--(47, 13.3125) -- (48, 13.3125)--(48, 12.375) -- (49, 12.375)--(49, 15.6875) -- (50, 15.6875)--(50, 18.0) -- (51, 18.0)--(51, 19.75) -- (52, 19.75)--(52, 17.5625) -- (53, 17.5625)--(53, 14.8125) -- (54, 14.8125)--(54, 15.9375) -- (55, 15.9375)--(55, 17.4375) -- (56, 17.4375)--(56, 14.25) -- (57, 14.25)--(57, 14.875) -- (58, 14.875)--(58, 18.4375) -- (59, 18.4375)--(59, 12.3125) -- (60, 12.3125)--(60, 14.6875) -- (61, 14.6875)--(61, 13.1875) -- (62, 13.1875)--(62, 14.375) -- (63, 14.375)--(63, 10.5625) -- (64, 10.5625)--(64, 7.875) -- (65, 7.875)--(65, 6.6875) -- (66, 6.6875)--(66, 12.0) -- (67, 12.0)--(67, 11.6875) -- (68, 11.6875)--(68, 4.8125) -- (69, 4.8125)--(69, 9.5) -- (70, 9.5)--(70, 13.9375) -- (71, 13.9375)--(71, 8.125) -- (72, 8.125)--(72, 13.5) -- (73, 13.5)--(73, 10.0) -- (74, 10.0)--(74, 17.0625) -- (75, 17.0625)--(75, 12.5625) -- (76, 12.5625)--(76, 19.0625) -- (77, 19.0625)--(77, 10.3125) -- (78, 10.3125)--(78, 14.1875) -- (79, 14.1875)--(79, 7.0625) -- (80, 7.0625)--(80, 17.0) -- (81, 17.0)--(81, 13.375) -- (82, 13.375)--(82, 13.125) -- (83, 13.125)--(83, 14.9375) -- (84, 14.9375)--(84, 15.75) -- (85, 15.75)--(85, 13.9375) -- (86, 13.9375)--(86, 10.4375) -- (87, 10.4375)--(87, 7.5) -- (88, 7.5)--(88, 14.375) -- (89, 14.375)--(89, 15.75) -- (90, 15.75)--(90, 18.3125) -- (91, 18.3125)--(91, 12.6875) -- (92, 12.6875)--(92, 11.125) -- (93, 11.125)--(93, 20.5625) -- (94, 20.5625)--(94, 22.8125) -- (95, 22.8125)--(95, 22.8125)
node[right] {};
    
\draw[thick, gray!80, dash pattern=on 28pt off 10pt, line width = 5pt]
(0,0) -- (28,0) -- (29,0) -- (29,0.8125) -- (30,0.8125) -- (30,1.625) -- (31,1.625) -- (31,2.4375) -- (32,2.4375) -- (32,3.25) -- (33,3.25) -- (33,4.0625) -- (34,4.0625) -- (34,4.875) -- (35,4.875) -- (35,5.6875) -- (36,5.6875) -- (36,6.5) -- (37,6.5) -- (37,7.3125) -- (38,7.3125) -- (38,8.125) -- (39,8.125) -- (39,8.9375) -- (40,8.9375) -- (40,9.75) -- (41,9.75) -- (41,10.5625) -- (42,10.5625) -- (42,11.375) -- (43,11.375) -- (43,12.1875) -- (44,12.1875) -- (44,13.0) -- (45,13.0) -- (45,13.8125) -- (46,13.8125) -- (46,14.625) -- (47,14.625) -- (47,15.4375) -- (48,15.4375) -- (48,16.25) -- (49,16.25) -- (49,17.0625) -- (50,17.0625) -- (50,16.25) -- (51,16.25) -- (51,15.4375) -- (52,15.4375) -- (52,14.625) -- (53,14.625) -- (53,13.8125) -- (54,13.8125) -- (54,13.0) -- (55,13.0) -- (55,12.1875) -- (56,12.1875) -- (56,11.375) -- (57,11.375) -- (57,10.5625) -- (58,10.5625) -- (58,9.75) -- (59,9.75) -- (59,8.9375) -- (60,8.9375) -- (60,8.125) -- (61,8.125) -- (61,7.3125) -- (62,7.3125) -- (62,6.5) -- (63,6.5) -- (63,5.6875) -- (64,5.6875) -- (64,4.875) -- (65,4.875) -- (65,4.0625) -- (66,4.0625) -- (66,3.25) -- (67,3.25) -- (67,2.5625) -- (68,2.5625) -- (68,1.625) -- (69,1.625) -- (69,0.8125) -- (70,0.8125) -- (70,0.0) -- (95,0.0)
node[right] {};

\draw[gray!50, dash pattern=on 70pt off 25pt, line width=5pt] (42.5,35) -- ++(5,0) node[right] {\scalebox{7}{\textbf{Energy budget}}};
\draw[black] (42.5,40) -- ++(5,0) node[right] {\scalebox{7}{\textbf{Consumption}}};
\draw[gray!80, dash pattern=on 28pt off 10pt, line width = 5pt] (42.5,30) -- ++(5,0) node[right] {\scalebox{7}{\textbf{Photovoltaic energy}}};

\end{tikzpicture}
        }
        \caption{Optimal schedule for variable consumptions.}
        \label{fig:6c}
    \end{subfigure}
    \caption{Comparison of solutions.}
    \label{fig:6}
\end{figure}
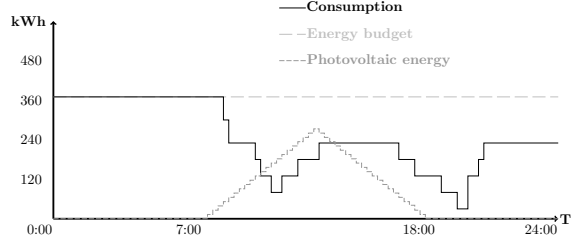
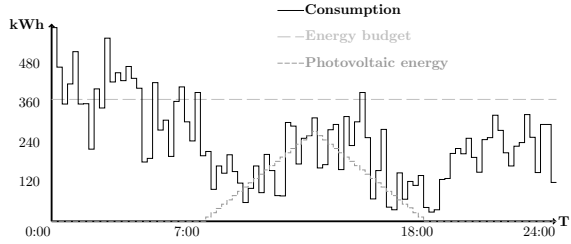
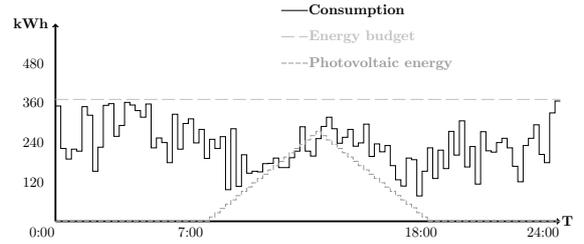

Figure~\ref{fig:6} illustrates this behavior by exploring in more detail a representative base configuration with the optimal schedules obtained for an instance with 20 jobs and 7 machines under 
fixed consumptions (Figure~\ref{fig:6a}), fixed consumptions but evaluated on (``plugged into'') variable consumptions (Figure~\ref{fig:6b}), and variable consumptions (Figure~\ref{fig:6c}).
In light gray, dark gray, and black we represent, respectively, the energy budget, the photovoltaic energy, and the total consumption in kWh of the schedule in each of the 96 time slots covering one day.
A direct comparison shows that the schedule obtained using fixed consumption tends to underestimate the peaks of energy usage, consuming more energy overnight at a lower cost, but exceeding the energy budget when evaluated on variable consumptions (Figure~\ref{fig:6b}). Consequently, the excess of energy purchased to complement the photovoltaic panels during the day (i.e., at higher prices) is reduced with respect to the energy obtained in the schedule with variable consumption functions, which in turn never exceeds the energy budget (Figure~\ref{fig:6c}). This type of effect is caused by neglecting the variability of energy consumption, resulting in schedules that may not be implemented in practice due to relevant practical limitations not considered in the models.
Incorporating time-dependent energy (or, more generally, resource) consumption functions in parallel machine scheduling models leads to more robust and cost-effective solutions.

\section{\rev{Conclusions}}
\label{sec:concl}

In this paper, we have studied the \problemname, an energy-efficient parallel machine scheduling problem inspired by a real-world application from the electric automotive industry. The \problemname\ introduces a novel element compared to the standard energy-efficient scheduling framework, as it accounts for variable energy consumption during job execution rather than assuming constant consumption. This variability, combined with a maximum energy budget constraint, significantly increases the complexity of the problem compared to classical parallel machine scheduling problems.

Our experimental framework revealed the critical impact of explicitly incorporating variable resource consumption within parallel machine scheduling.
The results indicate that replacing actual variable consumption with a constant average can lead to energy consumption exceeding the budget by up to 40\% on average. From a managerial perspective, ignoring variable consumption produces infeasible schedules with substantial energy budget violations, highlighting important practical implications.

From an algorithmic perspective, we presented an MILP model based on a time-indexed formulation and developed an Iterative Local Search matheuristic and a Simulated Annealing metaheuristic. The matheuristic integrates several constructive heuristics, local search operators, and an adapted MILP model. Computational experiments showed that while the MILP formulation performs well for small and medium instances, the matheuristic outperformed the MILP in solving large instances with variable consumption.
\rev{The matheuristic also outperformed the metaheuristic by finding better solutions, with statistically significant superiority on both the small/medium and large instances.}

This work opens several \rev{possible} future research directions. 
Methodologically, incorporating stochastic components to the model is a \rev{possible} direction. Photovoltaic energy production and, to some extent, job consumption are subject to uncertainty.
Scenario-based stochastic mathematical models and solution algorithms are an alternative to tackle these problems. Algorithmically, investigating other metaheuristic framework could lead to improved results. For instance, formulating subproblems in a rolling-horizon fashion similar to \cite{gaggero2023computational} is an interesting alternative.
From a modeling perspective, extending the \problemname\ to consider a continuous planning horizon with continuous variable consumption functions can provide deeper managerial insights and introduce new and challenging optimization problems into the energy-aware scheduling literature.

\vspace{-.5cm}
\section*{Acknowledgements}
This research was financially supported by Emilia-Romagna Region, Project DRIFT “DRive Into the FuTure – Studio e Ricerca di nuove tecnologie/metodologie operative, per una Battery Test House Intelligente, Innovativa e Sostenibile”. The academic work started when Juan José Miranda Bront was visiting the University of Modena and Reggio Emilia under the UNIMORE call Visiting Professor. We thank Reinova S.p.a. for supporting this research.
We acknowledge financial support from PNRR – M4C1 – Investimento 3.4, under project Green\&Pink for Sustainable Education - GPS Education Project Proposal TNE23-00012– CUP D74G23000280006.

\vspace{-.5cm}
\section*{Disclosure statement}
The authors also declare that they have no known competing financial interests or personal relationships that could have appeared to influence the work reported in this paper.

\vspace{-.5cm}
\section*{Data availability statement}
Data can be downloaded from the online repository: \href{https://github.com/regor-unimore/Energy-Efficient-Scheduling-Variable-Consumption.git}{https://github.com/regor-unimore/Energy-Efficient-Scheduling-Variable-Consumption.git}

\vspace{-.5cm}
\bibliographystyle{apalike}
\bibliography{itor}

\appendix
\vspace*{-10pt}
\section{Appendix A}\label{App:A}

Table~\ref{table:results-rectangular} shows the results produced by \milpmodel\ and \matheur\ for the small and medium instances with fixed consumption.
These instances are less challenging than the ones with variable consumptions, both in terms of feasibility and optimality (see Table \ref{table:results-random} for comparison). Both methods find a feasible solution for every instance. \milpmodel\ exhibits a remarkable behavior for this type of instances, reaching the optimal solution in 102 out of 108 instances, with execution times below 190 seconds on average. \matheur\ offers a viable alternative, finding near-optimal solutions with less than 0.1\% gap from the optimal, on average. We remark that adaptations of the MILP formulations by \cite{GAGGERO2023} that exploit the properties of fixed consumption functions may result in even shorter computation times for an exact approach.
\vspace{-0.5cm}
\begin{table}[ht]
\centering
\caption{Performance comparison of \milpmodel\ and \matheur\ on instances with fixed consumption. Values of \avglb\ and \avgub\ are expressed in cents of \euro. Values of \comptime\ and \inctime\ are expressed in seconds.}
\resizebox{\textwidth}{!}{%
\begin{tabular}{llc rrrrrr rrrrr}
\toprule
\multicolumn{3}{c}{Instances} & \multicolumn{6}{c}{\milpmodel} & \multicolumn{5}{c}{\matheur} \\
\cmidrule(lr){1-3} \cmidrule(lr){4-9} \cmidrule(lr){10-14}
$|T|$ & $(|J|$, $|I|)$ & \ninst & \nfeas & \nopt & \avglb & \avgub & \%gapLB & \comptime & \nfeas & \nopt & \%gapLB & \inctime & \%gapUB \\
\cmidrule(lr){1-3} \cmidrule(lr){4-9} \cmidrule(lr){10-14}
\multirow{3}*{48} 
& ([5, 20], [3, 7])  & 9 & 9 & 9 & 20560.83 & 20560.83 & 0.00 & 17.6 & 9 & 9 & 0.00 & 10.1 & 0.00 \\
& ([30, 90], [10, 20]) & 9 & 9 & 9 & 81020.83 & 81020.83 & 0.00 & 13.6 & 9 & 9 & 0.00 & 82.1 & 0.00 \\
& ([100, 200], [25, 35]) & 9 & 9 & 9 & 171448.33 & 171448.33 & 0.00 & 70.5 & 9 & 5 & 0.04 & 367.9 & 0.04 \\
\cmidrule(lr){2-3} \cmidrule(lr){4-9} \cmidrule(lr){10-14}
&  & 27 & 27 & 27 & 91010.00 & 91010.00 & 0.00 & 33.9 & 27 & 23 & 0.01 & 153.4 & 0.01 \\
\cmidrule(lr){1-3} \cmidrule(lr){4-9} \cmidrule(lr){10-14}
\multirow{3}*{72} 
& ([5, 20], [3, 7])  & 9 & 9 & 9 & 13049.04 & 13049.04 & 0.00 & 80.3 & 9 & 9 & 0.00 & 52.4 & 0.00 \\
& ([30, 90], [10, 20]) & 9 & 9 & 9 & 53728.15 & 53728.15 & 0.00 & 25.8 & 9 & 9 & 0.00 & 229.6 & 0.00 \\
& ([100, 200], [25, 35]) & 9 & 9 & 9 & 104639.33 & 104639.33 & 0.00 & 265.1 & 9 & 2 & 0.20 & 566.3 & 0.20 \\
\cmidrule(lr){2-3} \cmidrule(lr){4-9} \cmidrule(lr){10-14}
&  & 27 & 27 & 27 & 57138.84 & 57138.84 & 0.00 & 123.7 & 27 & 20 & 0.07 & 282.8 & 0.07 \\
\cmidrule(lr){1-3} \cmidrule(lr){4-9} \cmidrule(lr){10-14}
\multirow{3}*{96} 
& ([5, 20], [3, 7])  & 9 & 9 & 9 & 10033.96 & 10033.96 & 0.00 & 91.5 & 9 & 8 & 0.00 & 100.1 & 0.00 \\
& ([30, 90], [10, 20]) & 9 & 9 & 9 & 40698.50 & 40698.50 & 0.00 & 74.7 & 9 & 8 & 0.01  & 272.0 & 0.01 \\
& ([100, 200], [25, 35]) & 9 & 9 & 6 & 81255.23 & 81266.83 & 0.01 & 634.8 & 9 & 0 & 0.26 & 550.1 & 0.25\\
\cmidrule(lr){2-3} \cmidrule(lr){4-9} \cmidrule(lr){10-14}
&  & 27 & 27 & 24 & 43995.90 & 43999.76 & 0.00 & 267.0 & 27 & 16 &  0.09 & 307.4 &  0.09\\
\cmidrule(lr){1-3} \cmidrule(lr){4-9} \cmidrule(lr){10-14}
\multirow{3}*{120} 
& ([5, 20], [3, 7]) & 9 & 9 & 9 & 9642.27 & 9642.27 & 0.00 & 78.2 & 9 & 9 & 0.00 & 51.6 & 0.00 \\
& ([30, 90], [10, 20]) & 9 & 9 & 8 & 32813.41 & 32816.80 & 0.02 & 249.7 & 9 & 4 & 0.06  & 464.2 &  0.03\\
& ([100, 200], [25, 35]) & 9 & 9 & 7 & 70865.50 & 70870.40 & 0.01 & 632.9 & 9 & 0 & 0.35 & 548.9 & 0.35 \\
\cmidrule(lr){2-3} \cmidrule(lr){4-9} \cmidrule(lr){10-14}
&  & 27 & 27 & 24 & 37773.73 & 37776.49 & 0.01  & 320.3 & 27 & 13 & 0.14 & 354.9 &  0.13 \\
\cmidrule(lr){1-3} \cmidrule(lr){4-9} \cmidrule(lr){10-14}
Tot&  & 108 & 108 & 102 & 57479.62 & 57481.27 & 0.00 & 186.2 & 108 & 72 & 0.08 & 274.6 & 0.07 \\
\bottomrule
\end{tabular}%
}
\label{table:results-rectangular}
\end{table}

\section{\rev{Appendix B}}\label{App:B}
\rev{
Table~\ref{table:results-random} presents the complete results for \milpmodel, \matheur, and SA on small and medium instances (summarised in Table~\ref{table:results-random-summary} in the main text). The table reports 
the single execution of \milpmodel, and the average of ten executions of \matheur\ and SA with different random seeds, for the variable consumption instances in the corresponding row. For each method, \nfeas\ 
and \nopt\ indicate the number of instances for which it finds a feasible and an optimal solution, respectively, within the time limit.

For \milpmodel, we report the average values of the lower and upper bounds for each group of instances in \avglb\ and \avgub\ (expressed in cents of Euros, \euro), and the average final gap 
$\%\textrm{gapLB}$. \avglb, \avgub, and $\%\textrm{gapLB}$ are computed over a subset of instances within each group, whose size is indicated in 
parentheses next to each value; $\%\textrm{gapLB}$ is computed over the same subset as \avglb\ (i.e., when both a lower and an upper bound were found by \milpmodel).

For \matheur\ and SA, $\%\textrm{gapLB}$ is computed analogously with respect to the \milpmodel\ lower bound, as $100\times(z_{\textrm{ILS}} - z_{\textrm{lb}})/z_{\textrm{ILS}}$ and 
$100\times(z_{\textrm{SA}} - z_{\textrm{lb}})/z_{\textrm{SA}}$, respectively, over the same subset of instances used for \milpmodel's \avglb\ (in parentheses next to each value). \inctime\ indicates the time required to reach the best solution found during the execution; computation times for all three methods 
are reported in seconds. A direct comparison against \milpmodel's incumbent is provided in \%gapUB, calculated as $100\times(z_{\textrm{ILS}} - z_{\textrm{ub}})/z_{\textrm{ILS}}$ and $100\times(z_{\textrm{SA}} - z_{\textrm{ub}})/z_{\textrm{SA}}$, respectively, restricted to the subset of instances for 
which a feasible solution of \milpmodel\ exists (reported in column \#feas for \milpmodel, in parentheses next to each value). Gaps are computed on the same sets of instances across all executions of \matheur\ and SA.
}
\begin{sidewaystable}
\centering
\caption{Performance comparison of \milpmodel, \matheur, and SA on instances with variable consumption. 
Values of \avglb\ and \avgub\ are expressed in cents of \euro. Values of \comptime\ and \inctime\ are 
expressed in seconds. \rev{\matheur\ and SA values \nfeas, \nopt, \%gapLB, \inctime, and \%gapUB refer 
to the average results over 10 executions.} Numbers between parentheses for \avglb, \avgub, \%gapLB, and 
\%gapUB indicate the number of instances on which each value is computed in each execution.}
\resizebox{\textwidth}{!}{%
\begin{tabular}{llc rrrrrr rrrrr rrrrr}
\toprule
\multicolumn{3}{c}{Instances} & \multicolumn{6}{c}{\milpmodel} & \multicolumn{5}{c}{\matheur} & \multicolumn{5}{c}{\rev{SA}} \\
\cmidrule(lr){1-3} \cmidrule(lr){4-9} \cmidrule(lr){10-14} \cmidrule(lr){15-19}
$|T|$ & $(|J|$, $|I|)$ & \ninst & \nfeas & \nopt & \avglb & \avgub & \%gapLB & \comptime & \nfeas & \nopt & \%gapLB & \inctime & \%gapUB & \rev{\nfeas} & \rev{\nopt} & \rev{\%gapLB} & \rev{\inctime} & \rev{\%gapUB} \\
\cmidrule(lr){1-3} \cmidrule(lr){4-9} \cmidrule(lr){10-14} \cmidrule(lr){15-19}
\multirow{3}*{48} 
& ([5, 20], [3, 7]) & 9 & \textbf{7} & 4 & (7) 16891.80 & (7) 17441.14 & (7) 1.91 & 720.7 & \rev{\textbf{7.0}} & \rev{4.0} & \rev{(7) 1.88} & \rev{199.9} & \rev{(7) -0.03} & \rev{\textbf{7.0}} & \rev{3.0} & \rev{(7) 11.98} & \rev{8.3} & \rev{(7) 10.21} \\
& ([30, 90], [10, 20])  & 9 & \textbf{9} & 1 & (9) 79984.33 & (9) 80786.42 & (9) 1.06 & 1079.6 & \rev{\textbf{9.0}} & \rev{0.0} & \rev{(9) 1.24} & \rev{445.7} & \rev{(9) 0.19} & \rev{\textbf{9.0}} & \rev{0.0} & \rev{(9) 14.01} & \rev{170.2} & \rev{(9) 13.05} \\
& ([100, 200], [25, 35])  & 9 & \textbf{9} & 0 & (9) 168978.98 & (9) 169448.58 &  (9) 0.26 & 1200.0 & \rev{\textbf{9.0}} & \rev{0.0} & \rev{(9) 0.68} & \rev{522.9} & \rev{(9) 0.42} & \rev{\textbf{9.0}} & \rev{0.0} & \rev{(9) 11.91} & \rev{408.8} & \rev{(9) 11.67} \\
\cmidrule(lr){2-2} \cmidrule(lr){3-3} \cmidrule(lr){4-9} \cmidrule(lr){10-14} \cmidrule(lr){15-19}
&  & 27 & \textbf{25} & 5 & (25) 94356.49 & (25) 94968.12 & (25) 1.01 & 1000.1 &
\rev{\textbf{25.0}} & \rev{4.0} & \rev{(25) 1.22} & \rev{404.7} & \rev{(25) 0.21} & \rev{\textbf{25.0}} & \rev{3.0} & \rev{(25) 12.69} & \rev{210.8} & \rev{(25) 11.76} \\
\cmidrule(lr){1-3} \cmidrule(lr){4-9} \cmidrule(lr){10-14} \cmidrule(lr){15-19}
\multirow{3}*{72}
& ([5, 20], [3, 7]) & 9 & \textbf{8} & 3 & (8) 12258.37 & (8) 12675.12 & (8) 2.55 & 802.1 & \rev{\textbf{8.0}} & \rev{3.0} & \rev{(8) 2.56} & \rev{215.3} & \rev{(8) 0.01} & \rev{\textbf{8.0}} & \rev{1.5} & \rev{(8) 9.41} & \rev{44.8} & \rev{(8) 7.12} \\
& ([30, 90], [10, 20]) & 9 & \textbf{9} & 0 & (9) 53135.26 & (9) 53852.44 & (9) 1.17 & 1200.0 & \rev{\textbf{9.0}} & \rev{0.0} & \rev{(9) 1.48} & \rev{490.4} & \rev{(9) 0.31} & \rev{\textbf{9.0}} & \rev{0.0} & \rev{(9) 11.44} & \rev{124.9} & \rev{(9) 10.38} \\
& ([100, 200], [25, 35]) & 9 & 7 & 0 & (6) 90211.92 & (7) 95205.72 &  (6) 3.58 & 1200.0 & 
\rev{\textbf{9.0}} & \rev{0.0} & \rev{(6) 1.01} & \rev{570.8} & \rev{(7) -5.15} & \rev{\textbf{9.0}} & \rev{0.0} & \rev{(6) 16.36} & \rev{431.9} & \rev{(7) 11.19} \\
\cmidrule(lr){2-3} \cmidrule(lr){4-9} \cmidrule(lr){10-14} \cmidrule(lr){15-19}
&  & 27 & 24 & 3 & (23) 48589.39 & (24) 50317.71 & (23) 2.28 & 1067.4 & 
\rev{\textbf{26.0}} & \rev{3.0} & \rev{(23) 1.73} & \rev{433.6} & \rev{(24) -1.38} & \rev{\textbf{26.0}} & \rev{1.5} & \rev{(23) 12.02} & \rev{206.5} & \rev{(24) 9.53} \\
\cmidrule(lr){1-3} \cmidrule(lr){4-9} \cmidrule(lr){10-14} \cmidrule(lr){15-19}
\multirow{3}*{96} 
& ([5, 20], [3, 7]) & 9 & \textbf{6} & 3 & (6) 8790.69 & (6) 9102.25 & (6) 2.38 & 806.1 & 
\rev{\textbf{6.2}} & \rev{3.0} & \rev{(6) 2.59} & \rev{194.7} & \rev{(6) 0.22} & \rev{\textbf{6.2}} & \rev{1.1} & \rev{(6) 9.50} & \rev{60.1} & \rev{(6) 7.37} \\
& ([30, 90], [10, 20]) & 9 & 8 & 0 & (8) 36937.05 & (8) 37512.75 & (8) 1.54 & 1200.0 & 
\rev{\textbf{9.0}} & \rev{0.0} & \rev{(8) 1.75} & \rev{485.3} & \rev{(8) 0.21} & \rev{\textbf{9.0}} & \rev{0.0} & \rev{(8) 13.80} & \rev{180.6} & \rev{(8) 12.39} \\
& ([100, 200], [25, 35]) & 9 & 6 & 0 & (3) 66453.83 & (6) 73071.13 & (3) 7.57  & 1200.0 & 
\rev{\textbf{9.0}} & \rev{0.0} & \rev{(3) 1.56} & \rev{564.5} & \rev{(6) -16.90} & \rev{\textbf{9.0}} & \rev{0.0} & \rev{(3) 20.34} & \rev{475.0} & \rev{(6) 2.45} \\
\cmidrule(lr){2-3} \cmidrule(lr){4-9} \cmidrule(lr){10-14} \cmidrule(lr){15-19}
&  & 27 & 20 & 3 & (17) 32211.88 & (20) 33760.52 & (17) 2.90 &  1068.7 & 
\rev{\textbf{24.2}} & \rev{3.0} & \rev{(17) 2.01} & \rev{440.3} & \rev{(20) -4.92} & \rev{\textbf{24.2}} & \rev{1.1} & \rev{(17) 13.43} & \rev{259.2} & \rev{(20) 7.90} \\
\cmidrule(lr){1-3} \cmidrule(lr){4-9} \cmidrule(lr){10-14} \cmidrule(lr){15-19}
\multirow{3}*{120} 
& ([5, 20], [3, 7]) & 9 & \textbf{5} & 3 & (5) 6796.41 & (5) 6961.06 & (5) 1.71 & 848.6 & 
\rev{\textbf{5.0}} & \rev{3.0} & \rev{(5) 2.38} & \rev{246.9} & \rev{(5) 0.71} & \rev{\textbf{5.0}} & \rev{2.0} & \rev{(5) 8.39} & \rev{120.0} & \rev{(5) 6.85} \\
& ([30, 90], [10, 20]) & 9 & 5 & 0 & (5) 26239.79 & (5) 26572.90 & (5) 1.37  & 1200.0 & 
\rev{\textbf{6.7}} & \rev{0.0} & \rev{(5) 1.99} & \rev{516.5} & \rev{(5) 0.63} & \rev{\textbf{6.7}} & \rev{0.0} & \rev{(5) 21.67} & \rev{226.8} & \rev{(5) 20.52} \\
& ([100, 200], [25, 35]) & 9 & 3 & 0 & (1) 54839.46 & (3) 55288.44 &  (1) 0.81 & 1200.0 & 
\rev{\textbf{6.6}} & \rev{0.0} & \rev{(1) 1.02} & \rev{573.9} & \rev{(3) -16.33} & \rev{\textbf{6.5}} & \rev{0.0} & \rev{(1) 13.39} & \rev{479.4} & \rev{(3) 7.08} \\
\cmidrule(lr){2-3} \cmidrule(lr){4-9} \cmidrule(lr){10-14} \cmidrule(lr){15-19}
&  & 27 & 13 & 3 & (11) 20001.86 & (13) 20268.93 & (11) 1.47 & 1082.9 & 
\rev{\textbf{18.2}} & \rev{3.0} & \rev{(11) 2.08} & \rev{462.9} & \rev{(13) -3.25} & \rev{\textbf{18.2}} & \rev{2.0} & \rev{(11) 14.88} & \rev{287.7} & \rev{(13) 12.16} \\
\cmidrule(lr){1-3} \cmidrule(lr){4-9} \cmidrule(lr){10-14} \cmidrule(lr){15-19}
Tot & & 108 & 82 & 14 & (76) 55843.30 & (82) 56952.60 & (76) 1.88 & 1054.8 & 
\rev{\textbf{93.4}} & \rev{13.0} & \rev{(76) 1.68} & \rev{433.3} & \rev{(82) -2.06} & \rev{\textbf{93.4}} & \rev{7.6} & \rev{(76) 12.97} & \rev{237.1} & \rev{(82) 10.23} \\
\bottomrule
\end{tabular}%
}
\label{table:results-random}
\end{sidewaystable}

\section{\rev{Appendix C}}
\label{App C}

\rev{
Table~\ref{table:results-random-stats} reports, for the small and medium instances, the statistical analysis on the comparison of \matheur\ and SA.
The table reports RDPs of ILS and SA with respect to the best known solution (computed here across all three solution methods, including \milpmodel, reported in Table \ref{table:results-random} in Appendix \ref{App:B}), together with a Wilcoxon signed-rank test comparing the paired average objective values of \matheur\ and SA on each instance.
A Shapiro-Wilk test rejects normality of the paired per-instance differences ($p = 1.83\times10^{-7}$ overall, and $p \le 0.04$ for each individual value of $|T|$), so we apply the non-parametric Wilcoxon signed-rank test.
RPD is the relative percentage difference between the average objective value of a method over the 10 runs and the best known solution for that instance across all three solution methods (\milpmodel, SA, \matheur), i.e., $100\times(\bar{z}_{\text{method}} - z_{\text{BKS}})/z_{\text{BKS}}$. Numbers between parentheses for \%gapSA indicate the number of instances, out of \ninst, on which both SA and \matheur\ found a feasible solution in every one of the 10 runs, and on which \%gapSA and the significance test are therefore computed. Significance is based on a Wilcoxon signed-rank test on paired per-instance average objective values (\matheur\ vs.\ SA), two-sided: $^{*}p<0.05$, $^{**}p<0.01$, $^{***}p<0.001$.

Across the full set of 108 instances, \matheur\ achieves an average RPD of 0.38\%, against 14.42\% for SA, and the difference between the two methods is statistically significant ($p<0.001$) both overall and for each value of $|T|$, confirming that \matheur\ consistently and significantly outperforms SA also on small and medium instances.
}
 
\begin{table}[ht]
\centering
\caption{\rev{RPD and statistical significance analysis for SA and \matheur\ on the small and medium instances with variable energy consumption, aggregated by $|T|$. 
Values of Avg.\ Best are expressed in cents of \euro. Values of \inctime\ are expressed in seconds. }}
\setlength{\tabcolsep}{6pt}
\renewcommand{\arraystretch}{0.8}
\color{black}
\resizebox{\textwidth}{!}{%
\begin{tabular}{c r r rrrr rrrr}
\toprule
\multicolumn{2}{c}{\small Instances} & \multicolumn{1}{c}{\small \milpmodel} & \multicolumn{4}{c}{\small SA} & \multicolumn{4}{c}{\small \matheur} \\
\cmidrule(lr){1-2} \cmidrule(lr){3-3} \cmidrule(lr){4-7} \cmidrule(lr){8-11}
{\small $|T|$} & {\small \ninst} & {\small Avg.\ Best} & {\small \nfeas} & {\small Avg.\ Best} & {\small \rev{RPD (\%)}} & {\small \inctime} & {\small \nfeas} & {\small \inctime} & {\small \rev{RPD (\%)}} & {\small \%gapSA} \\
\cmidrule(lr){1-2} \cmidrule(lr){3-3} \cmidrule(lr){4-7} \cmidrule(lr){8-11}
{\small 48}  & {\small 27} & {\small 94968.12} & {\small 25} & {\small 106916.74} & {\small \rev{15.97}} & {\small 210.8} & {\small 25} & {\small 404.7} & {\small \rev{0.34}} & {\small (25) -12.22$^{***}$} \\
{\small 72}  & {\small 27} & {\small 54524.39} & {\small 26} & {\small 67127.70}  & {\small \rev{13.06}} & {\small 206.5} & {\small 26} & {\small 433.6} & {\small \rev{0.38}} & {\small (26) -14.58$^{***}$} \\
{\small 96}  & {\small 27} & {\small 44818.68} & {\small 24} & {\small 55452.30}  & {\small \rev{13.80}} & {\small 259.2} & {\small 24} & {\small 440.3} & {\small \rev{0.38}} & {\small (24) -15.01$^{***}$} \\
{\small 120} & {\small 27} & {\small 29102.07} & {\small 17} & {\small 41638.78}  & {\small \rev{15.04}} & {\small 287.7} & {\small 17} & {\small 462.9} & {\small \rev{0.45}} & {\small (17) -18.67$^{***}$} \\
\cmidrule(lr){1-2} \cmidrule(lr){3-3} \cmidrule(lr){4-7} \cmidrule(lr){8-11}
{\small Tot} & {\small 108} & {\small 60457.18} & {\small 92} & {\small 70184.27} & {\small \rev{14.42}} & {\small 237.1} & {\small 92} & {\small 433.3} & {\small \rev{0.38}} & {\small (92) -14.11$^{***}$} \\
\bottomrule
\end{tabular}%
}
\label{table:results-random-stats}
\end{table}

\end{document}